\documentclass[a4paper,12pt]{article}
\usepackage{longtable}
\title{Toward the classification\\of higher-dimensional toric Fano varieties}

\date{}

\renewcommand{\thefootnote}{\fnsymbol{footnote}}

\author{{\sc Hiroshi Sato}}

\newtheorem{Thm}{Theorem}[section]
\newtheorem{Prop}[Thm]{Proposition}
\newtheorem{Cor}[Thm]{Corollary}
\newtheorem{Conj}[Thm]{Conjecture}
\newtheorem{Lem}[Thm]{Lemma}

\newtheorem{Prob}[Thm]{Problem}
\newtheorem{Def}[Thm]{Definition}
\newtheorem{Rem}[Thm]{Remark}
\newtheorem{Ex}[Thm]{Example}

\newcommand{\proof}{Proof. \quad}
\newcommand{\qed}{\hfill q.e.d.}

\newcommand{\Hom}{\mathop{\rm Hom}\nolimits}

\newcommand{\G}{\mathop{\rm G}\nolimits}
\newcommand{\conv}{\mathop{\rm Conv}\nolimits}
\newcommand{\cone}{\mathop{\rm Cone}\nolimits}
\newcommand{\PC}{\mathop{\rm PC}\nolimits}
\newcommand{\relint}{\mathop{\rm Relint}\nolimits}
\newcommand{\orb}{\mathop{\rm orb}\nolimits}

\begin{document}
\maketitle

\renewcommand{\thefootnote}{}
\footnote{\hspace*{1.5em} {\em $1991$ Mathematics Subject Classification\/}.
Primary 14M25;
Secondary 14J45.}

\begin{abstract}                                                   

The purpose of this paper is to give basic tools for the classification of nonsingular toric Fano varieties by means of the notions of primitive collections and primitive relations due to Batyrev. By using them we can easily deal with equivariant blow-ups and blow-downs, and get an easy criterion to determine whether a given nonsingular toric variety is a Fano variety or not. As applications of these results, we get a toric version of a theorem of Mori, and can classify, in principle, all nonsingular toric Fano varieties obtained from a given nonsingular toric Fano variety by finite successions of equivariant blow-ups and blow-downs through nonsingular toric Fano varieties. Especially, we get a new method for the classification of nonsingular toric Fano varieties of dimension at most four. These methods are extended to the case of Gorenstein toric Fano varieties endowed with natural resolutions of singularities. Especially, we easily get a new method for the classification of Gorenstein toric Fano surfaces.

\end{abstract}

\section{Introduction}\label{intro}                                

\thispagestyle{empty}

\hspace{5mm} A {\em Gorenstein toric Fano variety} is a complete toric variety $X$ with at most Gorenstein singularities such that the anticanonical divisor $-K_X$ is ample. Gorenstein toric Fano varieties are very important as ambient spaces of Calabi-Yau varieties, and Batyrev \cite{batyrev2} systematically constructed examples of mirror symmetric pairs of Calabi-Yau varieties as hypersurfaces in Gorenstein toric Fano varieties. The set of isomorphism classes of Gorenstein toric Fano $d$-folds is a finite set for any dimension $d$ (see Batyrev \cite{batyrev1}). Nonsingular toric Fano $d$-folds are classified for $d\leq 4$ and Gorenstein toric Fano $d$-folds are classified for $d\leq 3$ (see Batyrev \cite{batyrev4} and Watanabe-Watanabe \cite{watanabe1} in the nonsingular cases, and Koelman \cite{koelman1}, Kreuzer-Skarke \cite{kreuzer1} and \cite{kreuzer2} in the Gorenstein cases). In this paper, we consider the classification of higher-dimensional nonsingular or Gorenstein toric Fano varieties using the notions of primitive collections and primitive relations introduced by Batyrev \cite{batyrev3}. First we consider the nonsingular case.

\begin{Def}\label{haita}

{\rm
Let $\mathcal{F}_d$ be the set of isomorphism classes of toric Fano $d$-folds. $X_1$ and $X_2$ in $\mathcal{F}_d$ are said to be {\em F-equivalent} if there exists a sequence of equivariant blow-ups and blow-downs from $X_1$ to $X_2$ through nonsingular toric {\em Fano} $d$-folds, namely there exist nonsingular toric Fano $d$-folds $Y_0=X_1,Y_1,\ldots,Y_{2l}=X_2$ together with finite successions $Y_j\rightarrow Y_{j-1}$ and $Y_j\rightarrow Y_{j+1}$, for each odd $1\leq j\leq 2l-1$, of equivariant blow-ups through nonsingular toric Fano $d$-folds. We denote the relation by $X_1\stackrel{\rm F}{\sim}X_2$. Then ``$\stackrel{\rm F}{\sim}$'' is obviously an equivalence relation.
}

\end{Def}

\begin{Rem}

{\rm
For equivariant birational maps of complete nonsingular toric varieties which need not be Fano varieties, related factorization conjectures have been proposed by Oda \cite{oda4}. The weak version analogous to the factorization in Definition \ref{haita} was proved by W\l odarczyk \cite{wlodarczyk1} and Morelli \cite{morelli1}, while the strong version was proved by Morelli \cite{morelli1} and later supplemented by Abramovich-Matsuki-Rashid \cite{abramovich1}.
}

\end{Rem}

As we see in this paper, if we get a complete system of representatives for $(\mathcal{F},\stackrel{\rm F}{\sim})$, then we get the classification of nonsingular toric Fano $d$-folds. The following conjecture for nonsingular toric Fano $d$-folds holds for $d\leq 4$ as a consequence of the known classification.

\begin{Conj}\label{mainconj}

Any nonsingular toric Fano $d$-fold is either pseudo-symmetric or F-equivalent to the $d$-dimensional projective space ${\bf P}^{d}$.

\end{Conj}

In this paper, we prove this conjecture for $d=3$ and $d=4$ without using the classification. As a result, we get a new method for the classification of nonsingular toric Fano $3$-folds and $4$-folds. Using this method for the classification, we can show that there exist 124 nonsingular toric Fano $4$-folds up to isomorphism. 

On the other hand, Gorenstein toric Fano $d$-folds are related to nonsingular toric weak Fano $d$-folds, where a nonsingular toric {\em weak} Fano variety is a nonsingular projective toric variety $X$ such that the anticanonical divisor $-K_X$ is nef and big, and the methods for nonsingular toric Fano $d$-folds are extended to the case of nonsingular weak toric Fano $d$-folds. As a result, we get a new method for the classification of Gorenstein toric Fano surfaces.

The content of this paper is as follows: In Section \ref{reflexive}, we study basic concepts on toric Fano varieties, and recall the correspondence between Gorenstein toric Fano varieties and reflexive polytopes. In Sections \ref{prim} and \ref{blow}, we introduce primitive collections and primitive relations. We can characterize toric Fano varieties using them, and calculate them before and after an equivariant blow-up. Moreover, we have a criterion for the possibility of an equivariant blow-down in terms of primitive collections and primitive relations. In Section \ref{blow}, we give a new nonsingular toric Fano $4$-fold which is missing in the classification of Batyrev \cite{batyrev4}. In Section \ref{decom}, we give a toric version of a theorem of Mori as an application of Sections \ref{prim} and \ref{blow}. In Section \ref{class}, we give a procedure for the classification which says that we have only to get a complete system of representatives for the F-equivalence relation for the set of isomorphism classes of nonsingular toric Fano $d$-folds. We also study a correspondence between toric weak Fano varieties and Gorenstein toric Fano varieties. Especially, we get a new method for the classification of Gorenstein toric Fano surfaces. In Sections \ref{nfano3} and \ref{nfano4}, we prove Conjecture \ref{mainconj} for $d=3$ and $d=4$. In Section \ref{four}, as an application of Sections \ref{prim} and \ref{blow}, we describe all the equivariant blow-up relations among nonsingular toric Fano $4$-folds using the classification of Batyrev \cite{batyrev4}.

The author wishes to thank Professors Tadao Oda, Yasuhiro Nakagawa and Takeshi Kajiwara for their advice and encouragement.

\section{Reflexive polytopes}\label{reflexive}                     

\hspace{5mm} In this section, we recall some basic notation and facts about toric Fano varieties (see Batyrev \cite{batyrev2}, Fulton \cite{fulton1}, and Oda \cite{oda2} for more details). The following notation is used throughout this paper.

Let $N$ be a free abelian group of rank $d$ and $M:=\Hom_{\bf Z}(N,{\bf Z})$ the dual group. The natural pairing $\langle\ ,\ \rangle :M\times N \rightarrow {\bf Z}$ is extended to a bilinear form $\langle\ ,\ \rangle :M_{\bf R}\times N_{\bf R} \rightarrow {\bf R}$ where $M_{\bf R}:=M\otimes_{\bf Z}{\bf R},\ N_{\bf R}:=N\otimes_{\bf Z}{\bf R}$.

For a finite complete fan $\Sigma$ in $N$ and $0\leq i\leq d$, we put $\Sigma(i):=\left\{ \sigma\in\Sigma\; |\;\dim\sigma=i\right\}$. Each $\tau\in\Sigma(1)$ determines a unique element $e(\tau)\in N$ which generates the semigroup $\tau\cap N$. We put
$$\G(\Sigma):=\left\{e(\tau)\in N\; |\;\tau\in\Sigma(1)\right\}$$
and $\G(\sigma):=\sigma\cap\G(\Sigma)$ for $\sigma\in \Sigma$.

\begin{Def}[Batyrev \cite{batyrev2}]

{\rm
A $d$-dimensional convex lattice polytope $\Delta\subset N_{\bf R}$ is called a {\em reflexive polytope} if the origin $0$ is in the interior of $\Delta$ and the polar
$$\Delta^{\ast}:=\left\{ y\in M_{\bf R}\; |\;\langle y,x\rangle \geq -1,\ \forall x\in \Delta\right\} \subset M_{\bf R}$$
is also a convex lattice polytope.
}

\end{Def}

For a $d$-dimensional convex polytope $\Delta\subset N_{\bf R}$ and $0\leq i\leq d-1$, we denote by $\Delta(i)$ the set of $i$-dimensional faces of $\Delta$.

Let $\Delta\subset N_{\bf R}$ be a convex lattice polytope such that $0$ is in the interior of $\Delta$. For any $i$-dimensional face $\delta\subset\Delta\ (0\leq i\leq d-1)$, let
$$\sigma(\delta):=\left\{ rx\in N_{\bf R}\; |\; r\in{\bf R}_{\geq 0},\ x\in\delta\right\}.$$
Then $\sigma(\delta)$ is an ($i+1$)-dimensional strongly convex rational polyhedral cone in $N_{\bf R}$. Moreover
$$\Sigma(\Delta):=\left\{\sigma(\delta)\; |\;\delta\in\Delta(i)\ (0\leq i\leq d-1)\right\}\cup\left\{ 0\right\}$$
is a finite complete fan in $N$.

\begin{Prop}[Batyrev \cite{batyrev2}]

If $\Delta\subset N_{\bf R}$ is a reflexive polytope, then $T_{N}{\rm emb}(\Sigma(\Delta))$ is a Gorenstein toric Fano variety. Conversely, if $\Sigma$ is a finite complete fan in $N$ such that $T_{N}{\rm emb}(\Sigma)$ is a Gorenstein toric Fano variety, then $\conv(\G(\Sigma))\subset N_{\bf R}$ is a reflexive polytope, where $\conv(\G(\Sigma))$ is the convex hull of $\G(\Sigma)\subset N_{\bf R}$. Moreover, any two Gorenstein toric Fano varieties $T_{N}{\rm emb}(\Sigma(\Delta_1))$ and $T_{N}{\rm emb}(\Sigma(\Delta_2))$ corresponding to two reflexive polytopes $\Delta_1\subset N_{\bf R}$ and $\Delta_2\subset N_{\bf R}$ are isomorphic if and only if $\Delta_1$ and $\Delta_2$ are equivalent up to unimodular transformation of the lattice $N$.

\end{Prop}

\begin{Rem}

{\rm
A reflexive polytope $\Delta$ is called a {\em Fano polytope} if $\Sigma(\Delta)$ is nonsingular.
}

\end{Rem}


\section{Primitive collections and primitive relations}\label{prim} 

\hspace{5mm} Primitive collections and primitive relations, introduced by Batyrev \cite{batyrev3}, are very convenient in describing higher-dimensional fans. So in this section, we recall these concepts and characterize toric Fano varieties using them.

\begin{Def}\label{pc1}

{\rm
Let $\Sigma$ be a finite complete {\em simplicial} fan in $N$. A nonempty subset $P\subset\G(\Sigma)$ is a {\em primitive collection} of $\Sigma$, if $\cone(P)\not\in\Sigma$, while $\cone\left(P\setminus\left\{ x\right\}\right)\in\Sigma$ for every $x\in P$, where $\cone(S):=\sum_{x\in S} {\bf R}_{\geq 0}x$ for any subset $S\subset N_{\bf R}$. 

We denote by $\PC(\Sigma)$ the set of primitive collections of $\Sigma$.
}

\end{Def}

\begin{Rem}

{\rm
By definition, for any subset $S\subset\G(\Sigma)$ which does not generate a cone in $\Sigma$, there exists a primitive collection $P\in\PC(\Sigma)$ such that $P\subset S$.
}

\end{Rem}

\begin{Def}\label{pc2}

{\rm
Let $\Sigma_1$ and $\Sigma_2$ be finite complete simplicial fans in $N$. Then $\PC(\Sigma_1)$ and $\PC(\Sigma_2)$ are {\em isomorphic} if there exists a bijective map $\varphi:\G(\Sigma_1)\rightarrow\G(\Sigma_2)$ which induces a well-defined bijective map
$$\varphi_{\ast}:\PC(\Sigma_1)\ni P\ \longmapsto\ \varphi(P)\in\PC(\Sigma_2).$$
}

\end{Def}

By Definitions \ref{pc1} and \ref{pc2}, we immediately get the following:

\begin{Prop}\label{combi}

Let $\Sigma_1$ and $\Sigma_2$ be finite complete simplicial fans in $N$. Then $\Sigma_1$ and $\Sigma_2$ are combinatorially equivalent if and only if $\PC(\Sigma_1)$ and $\PC(\Sigma_2)$ are isomorphic, where $\Sigma_1$ and $\Sigma_2$ are combinatorially equivalent if there exists a bijective map
$$\psi :\G(\Sigma_1)\ \longrightarrow\ \G(\Sigma_2)$$
such that for any nonempty subset $S\subset\G(\Sigma_1)$, we have $\cone(S)\in\Sigma_1$ if and only if $\cone(\psi(S))\in\Sigma_2$.

\end{Prop}

In the nonsingular case, we have the following additional information:

\begin{Def}\label{pr}

{\rm
Let $\Sigma$ be a finite complete nonsingular fan in $N$ and $P=\left\{ x_1,\ldots,x_l\right\}\in\PC(\Sigma)$. Then there is a unique element $\sigma(P)\in\Sigma$ such that
$$x_1+\cdots+x_l\in\relint(\sigma(P)),$$
where $\relint(S)$ is the relative interior of $S$ for any subset $S\subset N_{\bf R}$. 
Hence we get a linear relation
$$x_1+\cdots+x_l=a_1y_1+\cdots+a_my_m\ (a_1,\ldots,a_m\in{\bf Z}_{>0}),$$
where $\G(\sigma(P))=\left\{ y_1,\ldots,y_m\right\}$. We call this relation the {\em primitive relation} for $P$.

The integer $\deg P:=l-(a_1+\cdots+a_m)$ is called the {\em degree} of P.
}

\end{Def}

By this definition and Proposition \ref{combi}, we get the following characterization of isomorphism classes of complete nonsingular toric varieties.

\begin{Prop}\label{kero}

Let $\Sigma_1$ and $\Sigma_2$ be finite complete nonsingular fans in $N$. Then the complete nonsingular toric varieties $T_{N}{\rm emb}(\Sigma_1)$ and $T_{N}{\rm emb}(\Sigma_2)$ are isomorphic if and only if there exists an isomorphism from $\PC(\Sigma_1)$ to $\PC(\Sigma_2)$ which preserves their primitive relations.

\end{Prop}

Let $\Sigma$ be a finite complete nonsingular fan in $N$ and $X:=T_{N}{\rm emb}(\Sigma)$. Then for any $P\in\PC(\Sigma)$, we can define an element $r(P)\in A_1(X)$ in the following way, where $A_1(X)$ is the ${\bf Z}$-module of algebraic $1$-cycles modulo numerical equivalence.

\begin{Prop}[e.g., Fulton \cite{fulton1}, Oda \cite{oda2}]\label{exact}

Let $\Sigma$ be a finite complete nonsingular fan in $N$ and $X:=T_{N}{\rm emb}(\Sigma)$. Then we have an exact sequence of ${\bf Z}$-modules
$$0\ \longrightarrow \ M\ \stackrel{\varphi}{\longrightarrow} \ {\bf Z}^{\G(\Sigma)}\ \stackrel{\psi}{\longrightarrow} \ {\rm Pic}(X)\ \longrightarrow\ 0\ (\mbox{\rm exact}).$$

\end{Prop}

By the exact sequence in Proposition \ref{exact}, we have ${\rm Pic}(X)\cong{\bf Z}^{\G(\Sigma)}/M$ and hence
$$A_1(X)\cong\Hom_{\bf Z}({\rm Pic}(X),{\bf Z})\cong\Hom_{\bf Z}({\bf Z}^{\G(\Sigma)}/M,{\bf Z})\cong M^{\perp}\subset\Hom_{\bf Z}({\bf Z}^{\G(\Sigma)},{\bf Z}).$$
Consequently, we have
$$A_{1}(X)\cong \left\{(a_{x})_{x\in\G(\Sigma)}\in\Hom_{\bf Z}({\bf Z}^{\G(\Sigma)},{\bf Z})\;\left| \;\sum_{x\in \G(\Sigma)} a_{x}x=0\right.\right\}.$$
Let $P=\left\{x_1,\ldots,x_l\right\}\in\PC(\Sigma)$ and let
$$x_1+\cdots +x_l=a_1y_1+\cdots+a_my_m$$
be the primitive relation for $P$. Then we get a linear relation
$$x_1+\cdots+x_l-(a_1y_1+\cdots+a_my_m)=0.$$
Thus we can define $r(P)=(r(P)_x)_{x\in\G(\Sigma)}\in A_1(X)$ by
$$r(P)_x:=\left\{ \begin{array}{rcl}
1 & \mbox{if} & x=x_{i}\ (1\leq i \leq l) \\
-a_{j} & \mbox{if} & x=y_{j}\ (1\leq j\leq m) \\
0 & & \mbox{otherwise}. \\
\end{array} \right.$$

On the other hand, for any wall $\tau\in\Sigma(d-1)$, there is a linear relation
$$b_1z_1+\cdots+b_{d-1}z_{d-1}+b_dz_d+b_{d+1}z_{d+1}=0\ (b_1,\ldots,b_{d+1}\in{\bf Z},\ b_d=b_{d+1}=1),$$
where $\G(\tau)=\left\{z_1,\ldots,z_{d-1}\right\}$, while $\cone(\G(\tau)\cup\{z_d\})$ and $\cone(\G(\tau)\cup\{z_{d+1}\})$ are the $d$-dimensional strongly convex rational polyhedral cones in $\Sigma$ which contain $\tau$ as a face. We define $v(\tau)=(v(\tau)_x)_{x\in\G(\Sigma)}\in A_1(X)$ by
$$v(\tau)_x:=\left\{ \begin{array}{rcl}
b_{i} & \mbox{if} & x=z_{i}\ (1\leq i\leq d+1) \\
0 & & \mbox{otherwise}. \\
\end{array} \right. $$

Concerning this definition, the following is very useful.

\begin{Thm}[Batyrev \cite{batyrev3}, \cite{batyrev4}, Reid \cite{reid1}]\label{moricone}

Let $\Sigma$ be a finite complete nonsingular fan in $N$ and $X=T_{N}{\rm emb}(\Sigma)$. Then we have
$${\bf NE}(X)=\sum_{\tau\in\Sigma(d-1)}{\bf R}_{\geq 0}v(\tau)=\sum_{P\in\PC(\Sigma)}{\bf R}_{\geq 0}r(P)$$
where ${\bf NE}(X)\subset A_1(X)\otimes_{\bf Z}{\bf R}$ is the Mori cone of effective $1$-cycles.

\end{Thm}

The following theorem is the toric Nakai criterion.

\begin{Thm}[Oda \cite{oda2}, Oda-Park \cite{oda3}]\label{nakai}

Let $\Sigma$ be a finite complete nonsingular fan in $N$ and $X:=T_{N}{\rm emb}(\Sigma)$. Then a $T_{N}$-invariant divisor $D\in T_{N}{\rm Div}(X)$ is ample if and only if
$$\left( D.\overline{\orb(\tau)}\right) >0\mbox{ for all } \tau\in\Sigma(d-1).$$

\end{Thm}

By Theorems \ref{moricone} and \ref{nakai}, we can characterize nonsingular toric Fano varieties in terms of primitive collections.

\begin{Thm}[Batyrev \cite{batyrev4}]\label{fano}

Let $\Sigma$ be a finite complete nonsingular fan in $N$ and $X:=T_{N}{\rm emb}(\Sigma)$. Then $X$ is a nonsingular toric Fano variety $($resp. $-K_{X}$ is a nef divisor$)$ if and only if
$$\deg P>0\ (\mbox{resp.}\ \deg P\geq 0)\ \mbox{for\ all}\ P\in\PC(\Sigma).$$

\end{Thm}

\proof
${}^{t}(1,1,\ldots,1)\in{\bf Z}^{\G(\Sigma)}$ corresponds to the anticanonical divisor of $X$. So for $P\in\PC(\Sigma)$,
$$\left( -K_{X}.r(P)\right) =\deg P.$$
Hence by Theorems \ref{moricone} and \ref{nakai}, we are done.\hfill q.e.d.


\section{Equivariant blow-ups and blow-downs}\label{blow}           

\hspace{5mm} Let $\Sigma$ be a finite complete simplicial fan in $N$. In this section, we investigate how the set $\PC(\Sigma)$ of primitive collections change by star subdivisions. Especially we can deal with equivariant blow-ups and blow-downs of nonsingular complete toric varieties in terms of the primitive collections and primitive relations.

\begin{Def}\label{subdiv}

{\rm
Let $\Sigma$ be a finite complete simplicial fan in $N$ and $\sigma\in\Sigma$ with $\dim \sigma=l,\ 2\leq l\leq d$. For $x\in(\relint(\sigma))\cap N$ with $x$ primitive in $N$, we define the star subdivision of $\Sigma$ along $(\sigma,x)$ in the following way.

First, we define the strongly convex rational polyhedral cones $\sigma_i\ (1\leq i\leq l)$ by
$$\sigma_i:=\cone\left( \left\{x_1,\ldots,x_{i-1},x,x_{i+1},\ldots,x_l\right\} \right) \ (1\leq i\leq l)$$
where $\G(\sigma)=\left\{ x_1,\ldots,x_l\right\}$. Then for $\tau\in\Sigma$ such that $\sigma\prec\tau$, we can write $\tau$ uniquely as
$$\tau=\sigma+\tau'\mbox{ with }\tau'\in\Sigma,\ \sigma\cap\tau'=\{0\}.$$
In this notation, we have a finite complete simplicial fan $\Sigma^{\ast}_{(\sigma,x)}$ in $N$ defined by
$$\Sigma^{\ast}_{(\sigma,x)}:=\left( \Sigma\setminus\left\{ \tau\in\Sigma\; \left| \;\sigma\prec\tau\right. \right\} \right) \cup \left\{ \mbox{the faces of }\sigma_i+\tau'\; \left| \;\tau\in\Sigma,\ \sigma\prec\tau,\ 1\leq i\leq l\right. \right\}.$$We call $\Sigma^{\ast}_{(\sigma,x)}$ the star subdivision of $\Sigma$ along $(\sigma,x)$.
}

\end{Def}

\begin{Rem}[Fulton \cite{fulton1}, Oda \cite{oda2}]\label{oda}

{\rm
In {\rm Definition \ref{subdiv}}, if $\Sigma$ is nonsingular and $x=x_1+\cdots+x_l$, then the equivariant proper birational morphism $T_{N}{\rm emb}(\Sigma^{\ast}_{(\sigma,x)})\rightarrow T_{N}{\rm emb}(\Sigma)$ corresponding to this star subdivision is the equivariant blow-up along $\overline{\orb(\sigma)}$.
}

\end{Rem}

The following is one of the main theorems of this paper.

\begin{Thm}\label{primup}

Let $\Sigma$ be a finite complete simplicial fan in $N$, $\sigma\in\Sigma$ and $x$ a primitive element in $(\relint(\sigma))\cap N$. Then the primitive collections of $\Sigma^{\ast}_{(\sigma,x)}$ are

\begin{enumerate}

\item $\G(\sigma)$,
\item $P\in\PC(\Sigma)$ such that $\G(\sigma)\not\subset P$ and
\item the minimal elements in the set $\left\{ (P\setminus\G(\sigma))\cup\{x\}\; \left| \; P\in\PC(\Sigma),\ P\cap\G(\sigma)\neq\emptyset\right.\right\}$.

\end{enumerate}

\end{Thm}

To prove this theorem, we need the following three lemmas.

\begin{Lem}\label{easy}

Let $\Sigma$ be a finite complete simplicial fan in $N$, $\sigma\in\Sigma$ and $x$ a primitive element in $(\relint(\sigma))\cap N$. For any $\tau^{\ast}\in\Sigma^{\ast}_{(\sigma,x)}$, if $x\in\tau^{\ast}$ then there exists $\tau'\in\Sigma$ such that $\tau'\cap\sigma=\{0\}$ and $\tau^{\ast}\prec\sigma_i+\tau'\in\Sigma^{\ast}_{(\sigma,x)}$ for some $i$ with $1\leq i\leq l$. Moreover, $\sigma_j+\tau'\in\Sigma^{\ast}_{(\sigma,x)}$ for all $j$ $(1\leq j\leq l)$, where $l=\dim\sigma$.

\end{Lem}

The proof is trivial by Definition \ref{subdiv}.

\begin{Lem}\label{bunri}

Let $\Sigma$ be a finite complete simplicial fan in $N$, $\sigma\in\Sigma$ and $x$ a primitive element in $(\relint(\sigma))\cap N$. Then $P^{\ast}\in\PC\left( \Sigma^{\ast}_{(\sigma,x)}\right)$ and $x\in P^{\ast}$ imply $\G(\sigma)\cap P^{\ast}=\emptyset.$

\end{Lem}

\proof
Let $P^{\ast}\in\PC(\Sigma^{\ast}_{(\sigma,x)}),\ x\in P^{\ast}$ and suppose $\G(\sigma)\cap P^{\ast}\neq\emptyset$. Then $P^{\ast}\setminus\G(\sigma)$ generates a cone in $\Sigma$ containing $x$. So by Lemma \ref{easy}, there exists $\tau'\in\Sigma$ such that
$$P^{\ast}\setminus\G(\sigma)\subset\G(\sigma_i+\tau')\ (1\leq\exists i\leq l),\ \sigma\cap\tau'=\{0\},\ \sigma+\tau'\in\Sigma.$$
Since $(P^{\ast}\setminus\G(\sigma))\setminus\{x\}\subset\G(\tau')$, we have an index $j$ $(1\leq j\leq l)$ such that
$$P^{\ast}\subset\G(\sigma_j+\tau'),\ \sigma_j+\tau'\in\Sigma^{\ast}_{(\sigma,x)},$$
which contradicts the assumption.\hfill{\rm q.e.d.}

\begin{Lem}\label{sonzai}

Let $\Sigma$ be a finite complete simplicial fan in $N$, $\sigma\in\Sigma$ and $x$ a primitive element in $(\relint(\sigma))\cap N$. Then for any $P^{\ast}\in\PC(\Sigma^{\ast}_{(\sigma,x)})$ which contains $x$, there exists $P\in\PC(\Sigma)$ such that
$$\left( P\setminus\G(\sigma)\right) \cup\{x\}=P^{\ast}.$$

\end{Lem}

\proof
Let $P^{\ast}\in\PC(\Sigma^{\ast}_{(\sigma,x)}),\ x\in P^{\ast}$ and suppose $\G(\sigma)\cup\left( P^{\ast}\setminus\{x\}\right)$ generates a strongly convex rational polyhedral cone in $\Sigma$. Then there exists $\tau'\in\Sigma$ such that
$$\cone\left( \G(\sigma)\cup\left( P^{\ast}\setminus\{x\}\right) \right) =\sigma+\tau',\ \sigma\cap\tau'=\{0\}.$$
Since $\G(\sigma)\cap P^{\ast}=\emptyset$ by Lemma \ref{bunri}, we have $P^{\ast}\subset\G(\sigma_i+\tau')$ for all $i$ $(1\leq i\leq l)$. This contradicts $P^{\ast}\in\PC(\Sigma^{\ast}_{(\sigma,x)})$. Therefore $\G(\sigma)\cup\left( P^{\ast}\setminus\{x\}\right)$ contains a primitive collection of $\Sigma$.

Let $P\subset\G(\sigma)\cup\left( P^{\ast}\setminus\{x\}\right),\ P\in\PC(\Sigma)$. For any $y\in P^{\ast}\setminus\{x\}$, $P^{\ast}\setminus\{y\}$ generates a strongly convex rational polyhedral cone in $\Sigma^{\ast}_{(\sigma,x)}$ which contains $x$. Therefore by Lemma \ref{easy}, there exists $\tau'\in\Sigma$ such that
$$P^{\ast}\setminus\{y\}\subset\G(\sigma_i+\tau')\ (1\leq\exists i\leq l),\ \sigma\cap\tau'=\{0\}.$$
Then $P^{\ast}\setminus\{x,y\}\subset\G(\tau')$ because $\G(\sigma)\cap P^{\ast}=\emptyset$ by Lemma \ref{bunri}. So we have
$$\cone\left( \G(\sigma)\cup\left( P^{\ast}\setminus\{x,y\}\right) \right)=\sigma+\cone(P^{\ast}\setminus\{x,y\})\prec\sigma+\tau'\in\Sigma$$
and consequently $\G(\sigma)\cup\left( P^{\ast}\setminus\{x,y\}\right)$ generates a strongly convex rational polyhedral cone in $\Sigma$.

On the other hand, suppose $P^{\ast}\setminus\{x\}\not\subset P$. Then there exists $y\in P^{\ast}\setminus\{x\}$ such that $P\subset\G(\sigma)\cup\left( P^{\ast}\setminus\{x,y\}\right)$. This contradicts $P\in\PC(\Sigma)$. Therefore $P^{\ast}\setminus\{x\}\subset P$, hence clearly $\left( P\setminus\G(\sigma)\right) \cup\{x\}=P^{\ast}.$\hfill q.e.d.

\bigskip
We are now ready to prove Theorem \ref{primup}.

\bigskip
Proof of Theorem \ref{primup}. \quad
We put
$$\mathcal{P}:=\left\{ \left. P^{\ast}\in\PC(\Sigma^{\ast}_{(\sigma,x)})\; \right| \; x\not\in P^{\ast} \right\},\ \mathcal{P}':=\PC(\Sigma^{\ast}_{(\sigma,x)})\setminus\mathcal{P},$$
$$\mathcal{S}:=\left\{ P\in\PC(\Sigma)\; \left| \; \G(\sigma)\not\subset P\right. \right\} \cup\left\{ \G(\sigma)\right\}$$
and let $\mathcal{T}$ be the set of minimal elements of 
$$\left\{ \left( P\setminus\G(\sigma)\right) \cup\{x\}\; \left| \; P\in\PC(\Sigma),\ P\cap\G(\sigma)\neq\emptyset\right. \right\}.$$
Then to prove the theorem, we have only to prove $\mathcal{P}=\mathcal{S}$ and $\mathcal{P}'=\mathcal{T}$.

``$\mathcal{P}=\mathcal{S}$'' \quad
Trivially we have $\G(\sigma)\in\mathcal{P}$. Let $P\in\PC(\Sigma),\ \G(\sigma)\not\subset P$. Then for any $y\in P$, $P\setminus\{y\}$ generates a strongly convex rational polyhedral cone in $\Sigma^{\ast}_{(\sigma,x)}$ because $\G(\sigma)\not\subset P\setminus\{y\}$. On the other hand, since $x\not\in P$, $P$ does not generate a strongly convex rational polyhedral cone in $\Sigma^{\ast}_{(\sigma,x)}$. So we have $P\in\mathcal{P}$. Conversely, let $P^{\ast}\in\mathcal{P}$. If $\G(\sigma)\subset P^{\ast}$, then $P^{\ast}=\G(\sigma)\in\mathcal{S}$ since $\G(\sigma)\in\mathcal{P}$. If $\G(\sigma)\not\subset P^{\ast}$, then for any $y\in P^{\ast}$, $P^{\ast}\setminus\{y\}$ generates a strongly convex rational polyhedral cone in $\Sigma$ because $x\not\in P^{\ast}$. Clearly $P^{\ast}$ does not generate a strongly convex rational polyhedral cone in $\Sigma$. Therefore $P^{\ast}\in\PC(\Sigma)$ and we have $P^{\ast}\in\mathcal{S}$.

``$\mathcal{P}'=\mathcal{T}$'' \quad
Let $\left( P\setminus\G(\sigma)\right) \cup\{x\}\in\mathcal{T}\ \left( P\in\PC(\Sigma)\right)$ and suppose that $\left( P\setminus\G(\sigma)\right) \cup\{x\}$ generates a strongly convex rational polyhedral cone in $\Sigma^{\ast}_{(\sigma,x)}$. Then there exists $\tau'\in\Sigma$ such that
$$\cone\left( \left( P\setminus\G(\sigma)\right) \cup\{x\}\right) \prec\sigma_i+\tau'\in\Sigma^{\ast}_{(\sigma,x)}\ (1\leq\forall i\leq l),\ \sigma\cap\tau'=\{0\}.$$
Since $P\setminus\G(\sigma)\subset\G(\tau')$, we have $P\subset\G(\sigma+\tau')\ \left( \sigma+\tau'\in\Sigma\right)$, a contradiction to $P\in\PC(\Sigma)$. Therefore $\left( P\setminus\G(\sigma)\right) \cup\{x\}$ contains a primitive collection of $\Sigma^{\ast}_{(\sigma,x)}$. So let $P^{\ast}\subset\left( P\setminus\G(\sigma)\right) \cup\{x\},\ P^{\ast}\in\PC(\Sigma^{\ast}_{(\sigma,x)})$. Then $x\in P^{\ast}$ because $P\setminus\G(\sigma)$ generates a strongly convex rational polyhedral cone in $\Sigma^{\ast}_{(\sigma,x)}$. So by Lemma \ref{sonzai}, there exists $P'\in\PC(\Sigma)$ such that $P^{\ast}=\left( P'\setminus\G(\sigma)\right)\cup\{x\}$. Since $\left( P'\setminus\G(\sigma)\right)\cup\{x\}=P^{\ast}\subset\left( P\setminus\G(\sigma)\right)\cup\{x\}$, we have $\left( P'\setminus\G(\sigma)\right)\cup\{x\}=\left( P\setminus\G(\sigma)\right)\cup\{x\}$ by minimality. Therefore $\left( P\setminus\G(\sigma)\right)\cup\{x\}=P^{\ast}\in\PC(\Sigma^{\ast}_{(\sigma,x)})$. Conversely, let $P^{\ast}\in\PC(\Sigma^{\ast}_{(\sigma,x)}),\ x\in P^{\ast}$. Then by Lemma \ref{sonzai} $P^{\ast}$ is clearly expressed in the form as stated.\hfill q.e.d.

\bigskip

By using Theorem \ref{primup}, we can construct a nonsingular toric Fano $4$-fold which is missing in the table of Batyrev \cite{batyrev4}.

\begin{Ex}\label{newfano}

{\rm
Let $d=4$, $\Sigma$ a fan in $N$ corresponding to ${\bf P}^{2}\times{\bf P}^{2}$ and $\G(\Sigma)=\{x_1,\ldots,x_6\}$. Then the primitive relations of $\Sigma$ are$$x_1+x_2+x_3=0,\ x_4+x_5+x_6=0.$$
We get a nonsingular toric Fano $4$-fold $W$ by equivariant blow-ups of ${\bf P}^{2}\times{\bf P}^{2}$ along three $T_{N}$-invariant $2$-dimensional irreducible closed subvarieties
$$\overline{\orb\left( \left\{x_1,x_4\right\} \right) },\ \overline{\orb\left( \left\{x_2,x_5\right\} \right) },\ \overline{\orb\left( \left\{x_3,x_6\right\} \right) }.$$
Let $\Sigma_W$ be the fan in $N$ corresponding to $W$ and $\G(\Sigma_W)=\G(\Sigma)\cup\{x_7,x_8,x_9\}$. Then the primitive relations of $\Sigma_W$ are
$$x_1+x_4=x_7,\ x_2+x_5=x_8,\ x_3+x_6=x_9,$$
$$x_1+x_2+x_3=0,\ x_4+x_5+x_6=0,\ x_7+x_8+x_9=0,$$
$$x_1+x_2+x_9=x_6,\ x_4+x_5+x_9=x_3,\ x_1+x_3+x_8=x_5,$$
$$x_4+x_6+x_8=x_2,\ x_2+x_3+x_7=x_4,\ x_5+x_6+x_7=x_1,$$
$$x_1+x_8+x_9=x_5+x_6,\ x_4+x_8+x_9=x_2+x_3,\ x_2+x_7+x_9=x_4+x_6,$$
$$x_5+x_7+x_9=x_1+x_3,\ x_3+x_7+x_8=x_4+x_5,\ x_6+x_7+x_8=x_1+x_2.$$
This is easily confirmed by Theorem \ref{primup}. $W$ is missing in the table of Batyrev \cite{batyrev4}.

}

\end{Ex}

By Theorem \ref{primup}, we get a way to calculate $\PC(\Sigma^{\ast}_{(\sigma,x)})$ from $\PC(\Sigma)$. Conversely, by the following easy lemma, we get a way to calculate $\PC(\Sigma)$ from $\PC(\Sigma^{\ast}_{(\sigma,x)})$.

\begin{Lem}\label{forprimdown}

Let $\Sigma$ be a finite complete simplicial fan in $N$, $\sigma\in\Sigma$ and $x\in(\relint(\sigma))\cap N$ which generates the semigroup $({\bf R}_{\geq 0}x)\cap N$. If $P\in\PC(\Sigma)$ and $\G(\sigma)\subset P$, then $\left( P\setminus\G(\sigma)\right) \cup\{x\}\in\PC(\Sigma^{\ast}_{(\sigma,x)}).$

\end{Lem}

\proof
We have only to prove that $\left( P\setminus\G(\sigma)\right) \cup\{x\}$ is a minimal element in $\left\{ (P'\setminus\G(\sigma))\cup\{x\}\; \left| \; P'\in\PC(\Sigma),\ P'\cap\G(\sigma)\neq\emptyset\right.\right\}$. Suppose there exists $P'\in\PC(\Sigma)$ such that
$$P'\setminus\G(\sigma)\subset P\setminus\G(\sigma),\ P'\cap\G(\sigma)\neq\emptyset.$$
Since $\G(\sigma)\subset P$, we have $P'\subset P$, hence $P=P'$ because $P,\ P'\in\PC(\Sigma)$. Therefore $P$ is a minimal element.\hfill q.e.d.

\begin{Cor}\label{primdown}

Let $\Sigma$ be a finite complete simplicial fan in $N$, $\sigma\in\Sigma$ and $x\in(\relint(\sigma))\cap N$ which generates the semigroup $({\bf R}_{\geq 0}x)\cap N$. Then the primitive collections of $\Sigma$ are

\begin{enumerate}

\item $P^{\ast}\in\PC(\Sigma^{\ast}_{(\sigma,x)})$ such that $P^{\ast}\neq\G(\sigma),\ x\not\in P^{\ast}$ and
\item $\left( P^{\ast}\setminus\{x\}\right) \cup\G(\sigma)$ where $P^{\ast}\in\PC(\Sigma^{\ast}_{(\sigma,x)})$ such that $x\in P^{\ast}$ and $\left( P^{\ast}\setminus\{x\}\right) \cup S\not\in\PC(\Sigma^{\ast}_{(\sigma,x)})$ for any subset $S\subset\G(\sigma)$.

\end{enumerate}

\end{Cor}

This immediately follows from Theorem \ref{primup} and Lemma \ref{forprimdown}.

We end this section by giving an easy criterion for the possibility of equivariant blow-down in the nonsingular case.

\begin{Thm}\label{blowdown}

Let $\Sigma^{\ast}$ be a finite complete {\em nonsingular} fan in $N$. Then the following are equivalent.

\begin{enumerate}

\item There exist a complete nonsingular toric variety $X$ and an equivariant blow-up $\varphi:T_{N}{\rm emb}(\Sigma^{\ast})\rightarrow X$ along a $T_{N}$-invariant closed irreducible subvariety of $X$.
\item There exists $P^{\ast}\in\PC(\Sigma^{\ast})$ such that the corresponding primitive relation is
$$x_1+\cdots+x_l=x,\ P^{\ast}=\{x_1,\ldots,x_l\},\mbox{ for some }x\in\G(\Sigma^{\ast})$$
and for any $\sigma^{\ast}\in\Sigma^{\ast}$ which contains $x$, each of
$$\left( \G(\sigma^{\ast})\cup P^{\ast}\right) \setminus\{x_i\}\qquad (1\leq\forall i\leq l)$$
generates a strongly convex rational polyhedral cone in $\Sigma^{\ast}$.
\item There exists $P^{\ast}\in\PC(\Sigma^{\ast})$ such that the corresponding primitive relation is
$$x_1+\cdots+x_l=x,\ P^{\ast}=\{x_1,\ldots,x_l\},\mbox{ for some }x\in\G(\Sigma^{\ast})$$
and for any $P'\in\PC(\Sigma^{\ast})$ which satisfies the conditions $P^{\ast}\cap P'\neq\emptyset$ and $P^{\ast}\neq P'$,
$$(P'\setminus P^{\ast})\cup\{x\}$$
contains a primitive collection of $\Sigma^{\ast}$.

\end{enumerate}

\end{Thm}

\proof
We prove $(1)\Longrightarrow(3)\Longrightarrow(2)\Longrightarrow(1)$.

$(1)\Longrightarrow(3)$ is trivial by Theorem \ref{primup}.

$(3)\Longrightarrow(2).$ Suppose that there exists $\sigma^{\ast}\in\Sigma^{\ast}$ such that $x\in\sigma^{\ast}$ and
$$\left( \G(\sigma^{\ast})\cup P^{\ast}\right) \setminus\{x_i\}\mbox{ for some }i\ (1\leq i\leq l)$$
does not generate a strongly convex rational polyhedral cone in $\Sigma^{\ast}$. Then $\left( \G(\sigma^{\ast})\cup P^{\ast}\right) \setminus\{x_i\}$ contains a primitive collection $P'\in\PC(\Sigma^{\ast})$. Since $P^{\ast}\cap P'\neq\emptyset$ and $P^{\ast}\neq P'$, by $(3)$,
$$\left( P'\setminus P^{\ast}\right) \cup\{x\}\subset\G(\sigma^{\ast})$$
contains a primitive collection of $\Sigma^{\ast}$, a contradiction.

$(2)\Longrightarrow(1)$. For any $\sigma^{\ast}\in\Sigma^{\ast}$ which contains $x$, define a strongly convex rational polyhedral cone $\sigma'$ in $N_{\bf R}$ by
$$\sigma':=\cone\left( \left( \G(\sigma^{\ast})\cup P^{\ast}\right) \setminus\{x\}\right).$$
Then the finite complete nonsingular fan $\Sigma$ in $N$ defined by
$$\Sigma:=\left( \Sigma^{\ast}\setminus\left\{ \sigma^{\ast}\in\Sigma^{\ast}\left| x\in\sigma^{\ast}\right. \right\} \right) \cup\left\{ \sigma'\mbox{ and the faces of }\sigma'\left| \sigma^{\ast}\in\Sigma^{\ast},\ x\in\sigma^{\ast}\right. \right\}$$
gives a complete nonsingular toric variety $X=T_{N}{\rm emb}(\Sigma)$ and an equivariant blow-up $\varphi:T_{N}{\rm emb}(\Sigma)\rightarrow X$.\hfill q.e.d.

\bigskip
The equivalence $(1)\Longleftrightarrow(3)$ is a useful criterion for the possibility of equivariant blow-down in the nonsingular case.


\section{Decomposition of birational morphisms}\label{decom}        

\hspace{5mm} In this section, we prove a toric version of a theorem of Mori which says ``a proper birational morphism between nonsingular Fano $3$-folds is always decomposed into a composite of blow-ups'', and consider the higher-dimensional version. In the proof of this theorem, the results of Sections \ref{prim} and \ref{blow} are used.

The following proposition is important in proving the main theorem in this section.

\begin{Prop}\label{keylemma}

Let $X:=T_N{\rm emb}(\Sigma)$ be a nonsingular toric Fano $d$-fold $($resp. $-K_X$ is nef$)$, $x_1+\cdots+x_l=x$ a primitive relation of $\Sigma$ and $\varphi:X\rightarrow X':=T_{N}{\rm emb}(\Sigma')$ the equivariant blow-down with respect to $x_1+\cdots+x_l=x$. Then $X'$ is not a nonsingular toric Fano $d$-fold $($resp. $-K_{X'}$ is not nef$)$ if and only if there exists a primitive relation of $\Sigma$ of the form
$$y_1+\cdots+y_m=a_1z_1+\cdots+a_nz_n+bx+c_1x_1+\cdots+c_{l-1}x_{l-1}$$
up to change of the indices, such that
\begin{enumerate}
\item $a_1,\ldots,a_n,b>0,\ c_1,\ldots,c_{l-1}\geq 0$,
\item $m-(a_1+\cdots+a_n+b+c_1+\cdots+c_{l-1})>0$ $($resp. $\geq 0)$,
\item $m-(a_1+\cdots+a_n+bl+c_1+\cdots+c_{l-1})\leq 0$ $($resp. $<0)$ and
\item $m+n+l\leq d+1$.
\end{enumerate}

\end{Prop}

\proof
The sufficiency is trivial by Theorem \ref{fano}.

By Corollary \ref{primdown}, for any new primitive collection $P'\in\PC(\Sigma')$ added by the equivariant blow-down with respect to $x_1+\cdots+x_l=x$, there exists
$$P=\{u_1,\ldots,u_r,x\}\in\PC(\Sigma)$$
such that $P'=\{u_1,\ldots,u_r,x_1,\ldots,x_l\}$. Let the primitive relation corresponding to $P$ be
$$u_1+\cdots+u_r+x=h_1v_1+\cdots+h_sv_s.$$
Then $\cone\left( \{v_1,\ldots,v_s\}\right)\in\Sigma'$ because $x\not\in\{v_1,\ldots,v_s\}$. So the primitive relation corresponding to $P'$ is
$$u_1+\cdots+u_r+x_1+\cdots+x_l=h_1v_1+\cdots+h_sv_s.$$
Therefore $\deg P'=r+l-(h_1+\cdots+h_s)>r+1-(h_1+\cdots+h_s)=\deg P>0$.

By the above discussion, if $X'$ is not a nonsingular toric Fano $d$-fold, then there exists a primitive collection $P$ in $\PC(\Sigma)$ such that $P$ is in $\PC(\Sigma')$, its primitive relation contains $x$ on the right-hand side and $r(P)$ is contained in an extremal ray of ${\bf NE}(X')$. So we get the conditions $(1)$ and $(4)$. Because $X$ is a Fano variety while $X'$ is not a Fano variety, we get the conditions $(2)$ and $(3)$.\hfill q.e.d.

\bigskip

\begin{Ex}\label{youming}

{\rm

We consider Proposition \ref{keylemma} in the case of the equivariant blow-down $\varphi:X\rightarrow X'$ with respect to the primitive relation of $\Sigma$ of the form $x_1+x_2=x.$

\begin{enumerate}

\item ``$d=2$'' \quad
$X'$ is always a nonsingular toric Fano surface. On the other hand, if $-K_X$ is nef, then $-K_{X'}$ is always nef.
\item ``$d=3$'' \quad
$X'$ is not a nonsingular toric Fano $3$-fold if and only if there exists the following primitive relation of $\Sigma$.
$$y_1+y_2=x\quad \left( \left\{ y_1,y_2\right\} \cap \left\{ x_1,x_2\right\} =\emptyset\right) .$$
\item ``$d=4$'' \quad
$X'$ is not a nonsingular toric Fano $4$-fold if and only if there exists one of the following primitive relations of $\Sigma$.
$$y_1+y_2=x,\ y_1+y_2+y_3=2x,\ y_1+y_2+y_3=x+x_1\quad \left( \left\{ y_1,y_2,y_3\right\} \cap \left\{ x_1,x_2\right\} =\emptyset\right) .$$

\end{enumerate}

Next let $d=3$ and let $\varphi:X\rightarrow X'$ be the equivariant blow-down with respect to the primitive relation of $\Sigma$, $x_1+x_2+x_3=x.$ Then $X'$ is always a nonsingular toric Fano $3$-fold by Proposition \ref{keylemma}.

We need these facts later.

}

\end{Ex}

\bigskip

The following is the toric version of the Mori theory.

\begin{Prop}[Reid \cite{reid1}]\label{contraction}

Let $\Sigma$ be a finite complete nonsingular fan in $N$, $X:=T_{N}{\rm emb}(\Sigma)$ a projective toric variety, and $P=\{x_1,\ldots,x_l\}\in\PC(\Sigma)$ with the primitive relation corresponding to $P$ being $x_1+\cdots+x_l=a_1y_1+\cdots+a_my_m$. If $r(P)$ is contained in an extremal ray of ${\bf NE}(X)$ and $m\geq 1$, then there exist a nonsingular projective toric $d$-fold $X'$ and an equivariant morphism
$${\rm Cont}_{P}:X\longrightarrow X'$$
such that the following are satisfied:

\begin{enumerate}

\item For any $\tau\in\Sigma$, the image of $\overline{{\rm orb}(\tau)}$ by ${\rm Cont}_{P}$ is a point if and only if $v(\tau)=r(P)\in A_1(X).$
\item Let $\Sigma'$ be a fan in $N$ such that $X'=T_{N}{\rm emb}(\Sigma')$. If $m=1$ then $\Sigma'$ is simplicial and
$$\sigma'=\cone\left( \{x_1,\ldots,x_l\}\right) \in\Sigma',\ \G(\Sigma')=\G(\Sigma)\setminus\{y_1\}.$$
Moreover $\Sigma=(\Sigma')^{\ast}_{(\sigma',x)}$ where $x:=(x_1+\cdots+x_l)/a_1$. Especially if $a_1=1$, then $X'$ is nonsingular and ${\rm Cont}_{P}$ is an equivariant blow-up.

\end{enumerate} 

\end{Prop}

To prove the main theorem in this section, we suppose $d=3$. Let $\varphi:Y\longrightarrow X$ be an equivariant morphism between nonsingular toric Fano $3$-folds, and $\Sigma$ and $\widetilde{\Sigma}$ fans in $N$ such that $X=T_{N}{\rm emb}(\Sigma)$ and $Y=T_{N}{\rm emb}(\widetilde{\Sigma})$. To apply Propositions \ref{keylemma} and \ref{contraction}, we have to investigate the subdivision of a $3$-dimensional strongly convex rational polyhedral cone in $\Sigma$. The following lemma is fundamental in classifying subdivisions.

\begin{Lem}\label{piyo}

Let $d={\rm rank}N=3$, $\Sigma$ and $\widetilde{\Sigma}$ finite complete nonsingular fans in $N$ and $\varphi:T_{N}{\rm emb}(\widetilde{\Sigma}) \longrightarrow T_{N}{\rm emb}(\Sigma)$ an equivariant morphism. For any $\sigma\in\Sigma(3)$ such that $\G(\sigma)=\{x_1,x_2,x_3\}$, let $\widetilde{\sigma}$ be the unique strongly convex rational polyhedral cone in $\widetilde{\Sigma}\setminus\{0\}$ such that $x_1+x_2+x_3\in\relint(\widetilde{\sigma})$. Then we have the following.

\begin{enumerate}

\item $\dim\widetilde{\sigma}=3\Longleftrightarrow \sigma=\widetilde{\sigma}\in\widetilde{\Sigma}.$
\item $\dim\widetilde{\sigma}=2\Longleftrightarrow \G(\widetilde{\sigma})=\{x,x_3\}$ where $x:=x_1+x_2$ up to change of the indices.
\item $\dim\widetilde{\sigma}=1\Longleftrightarrow \G(\widetilde{\sigma})=\{x\}$ where $x:=x_1+x_2+x_3$.

\end{enumerate}

\end{Lem}

\proof
The sufficiency is trivial. Let $s=\dim\widetilde{\sigma}$ and $\G(\widetilde{\sigma})=\{y_1,\ldots,y_s\}$. Then $\widetilde{\sigma}\subset\sigma$ since $\varphi$ is an equivariant morphism, and so we have
$$y_i=a_{i1}x_1+a_{i2}x_2+a_{i3}x_3\ (1\leq i\leq s),\ a_{ij}\in{\bf Z}_{\geq 0}\ (1\leq i\leq s,\ 1\leq j\leq 3).$$
If we put $x_1+x_2+x_3=b_1y_1+\cdots+b_sy_s\ (b_1,\ldots,b_s\in{\bf Z}_{>0})$, then $b_1=\cdots=b_s=1$ because $a_{ij}\ (1\leq i\leq s,\ 1\leq j\leq 3)$ are nonnegative.\hfill q.e.d.

\bigskip
Now we are ready to classify the subdivisions of a $3$-dimensional strongly convex rational polyhedral cone $\sigma\in\Sigma(3)$. There are five types of subdivisions for $\sigma$. Let $\G(\sigma)=\{x_1,x_2,x_3\}$.

\begin{enumerate}

\item ``$\dim\widetilde{\sigma}=3$'' \quad
$\sigma=\widetilde{\sigma}\in\widetilde{\Sigma}(3)$ by Lemma \ref{piyo}.
\item ``$\dim\widetilde{\sigma}=2$'' \quad
By Lemma \ref{piyo}, we have $x_1+x_2+x_3\in\cone\left( \{x_3,x_4\}\right) \in\widetilde{\Sigma}(2)$, where $x_4:=x_1+x_2\in\G(\widetilde{\Sigma})$. Then $\{x_1,x_2\}\in\PC(\widetilde{\Sigma})$ and $r(\{x_1,x_2\})$ is contained in an extremal ray of ${\bf NE}(Y)$ since $\deg\left( \{x_1,x_2\}\right) =1$. So $\sigma_1:=\cone\left( \{x_1,x_3,x_4\}\right) ,\ \sigma_2:=\cone\left( \{x_2,x_3,x_4\}\right)$ are in $\widetilde{\Sigma}(3)$ and $\sigma=\sigma_1\cup\sigma_2$ by Theorem \ref{blowdown} $(2)$ and Proposition \ref{contraction}.
\item ``$\dim\widetilde{\sigma}=1$ and $\{x_1,x_2\}\in\PC(\widetilde{\Sigma})$'' \quad
Let $x_4:=x_1+x_2+x_3\in\G(\widetilde{\Sigma})$ and $x_5:=x_1+x_2$. Then $x_5\in\G(\widetilde{\Sigma})$ and the primitive relation corresponding to $\{x_1,x_2\}$ is $x_1+x_2=x_5$. Since $x_3+x_5=x_4$, we have $\{x_3,x_5\}\in\PC(\widetilde{\Sigma})$ and $x_3+x_5=x_4$ is the corresponding primitive relation. So since $r(\{x_1,x_2\})$ and $r(\{x_3,x_5\})$ are contained in an extremal ray of ${\bf NE}(Y)$, we see that $\sigma_1:=\cone\left( \{x_1,x_3,x_4\}\right) ,\ \sigma_2:=\cone\left( \{x_1,x_4,x_5\}\right) ,\ \sigma_3:=\cone\left( \{x_2,x_3,x_4\}\right) ,\ \sigma_4:=\cone\left( \{x_2,x_4,x_5\}\right)$ are in $\widetilde{\Sigma}(3)$ and $\sigma=\sigma_1\cup\sigma_2\cup\sigma_3\cup\sigma_4$ for the same reason as above.
\item ``$\dim\widetilde{\sigma}=1$, $\{x_1,x_2,x_3\}\in\PC(\widetilde{\Sigma})$ and $r\left( \{x_1,x_2,x_3\}\right)$ is contained in an extremal ray of ${\bf NE}(Y)$'' \quad
Let $x_4:=x_1+x_2+x_3$. Then by Proposition \ref{contraction}, we have $\sigma_1:=\cone\left( \{x_1,x_2,x_4\}\right) ,\ \sigma_2:=\cone\left( \{x_2,x_3,x_4\}\right) ,\ \sigma_3:=\cone\left( \{x_1,x_3,x_4\}\right)$ are in $\widetilde{\Sigma}(3)$ and $\sigma=\sigma_1\cup\sigma_2\cup\sigma_3$.
\item ``$\dim\widetilde{\sigma}=1$, $\{x_1,x_2,x_3\}\in\PC(\widetilde{\Sigma})$ and $r\left( \{x_1,x_2,x_3\}\right)$ is not contained in an extremal ray of ${\bf NE}(Y)$'' \quad
Let $x_4:=x_1+x_2+x_3$. Then the primitive relation corresponding to $\{x_1,x_2,x_3\}$ is $x_1+x_2+x_3=x_4$ and so $\deg \left( \{x_1,x_2,x_3\}\right) =2$. Therefore there exist two primitive collections $P_1,\ P_2\in\PC(\widetilde{\Sigma})$ such that $\deg P_1=\deg P_2=1$ and $r\left( \{x_1,x_2,x_3\}\right) =r(P_1)+r(P_2)$. On the other hand, there are two types of primitive relations corresponding to the primitive collection $P$ such that $\deg P=1$ and $r(P)$ is contained in an extremal ray. The possibilities are
$$(a)\ z_1+z_2+z_3=2z_4,\ (b)\ w_1+w_2=w_3.$$
By easy calculation, the combinations $\left( (a),(a)\right)$ and $\left( (b),(b)\right)$ are impossible. In the case of the combination $\left( (a),(b)\right)$, we have $z_4=w_1=x_4,\ w_3=z_1,\ w_2=x_1,\ z_2=x_2$ and $z_3=x_3$. Then putting $x_5:=z_1$, we have $\sigma_1:=\cone\left( \{x_1,x_2,x_5\}\right) ,\ \sigma_2:=\cone\left( \{x_2,x_4,x_5\}\right) ,\ \sigma_3:=\cone\left( \{x_1,x_3,x_5\}\right) ,\ \sigma_4:=\cone\left( \{x_3,x_4,x_5\}\right) ,\ \sigma_5:=\cone\left( \{x_2,x_3,x_4\}\right)$ are in $\widetilde{\Sigma}(3)$ and $\sigma=\sigma_1\cup\sigma_2\cup\sigma_3\cup\sigma_4\cup\sigma_5$ for the same reason as in $(2)$.

\end{enumerate}

By the above classification, we get the following main theorem in this section. This is a toric version of a theorem of Mori.

\begin{Thm}\label{pikachu}

Let $X$ and $Y$ be nonsingular toric Fano $3$-folds, and $\varphi:Y\longrightarrow X$ an equivariant morphism. Then we have a decomposition of $\varphi$
$$Y=X_r\stackrel{\varphi_r}{\longrightarrow}X_{r-1}\stackrel{\varphi_{r-1}}{\longrightarrow}\cdots\cdots\stackrel{\varphi_2}{\longrightarrow}X_1\stackrel{\varphi_1}{\longrightarrow}X_0=X$$
where $X_i\ (0\leq i\leq r)$ is a nonsingular toric Fano $3$-fold, $\varphi_j\ (2\leq j\leq r)$ is an equivariant blow-up along a $T_{N}$-invariant $1$-dimensional irreducible closed subvariety of $X_{j-1}$ and $\varphi_1$ is an equivariant blow-up along some $T_{N}$-invariant points of $X$.

\end{Thm}

\proof
In the above classification, carry out equivariant blow-downs in the order $(3)\Longrightarrow(2)\Longrightarrow(1)$, $(2)\Longrightarrow(1)$, $(5)\Longrightarrow(4)\Longrightarrow(1)$ and $(4)\Longrightarrow(1)$. Then by Proposition \ref{keylemma} and Example \ref{youming}, we get a decomposition as in the statement.\hfill q.e.d.

\bigskip
If $d\geq 4$, the method we employed in the $3$-dimensional case is insufficient. For example, in the case of $d=4$, there is a subdivision of a $4$-dimensional strongly convex rational polyhedral cone $\sigma\in\Sigma(4)$ such that the primitive relations corresponding to $\left\{ P\in\PC(\widetilde{\Sigma})\; \left| \; P\subset\sigma\right. \right\}\subset\PC(\widetilde{\Sigma})$ are
$$x_1+x_2+x_3=x_5,\ x_2+x_4=x_6\mbox{ and }x_1+x_3+x_6=x_4+x_5,$$
where $\G(\sigma)=\{x_1,x_2,x_3,x_4\},\ x_5,\ x_6\in\G(\widetilde{\Sigma})$. This does not contradict the fact that $Y$ is a nonsingular toric Fano variety, but we cannot tell by Proposition \ref{keylemma} whether the equivariant blow-down of $Y$ with respect to the primitive relation $x_2+x_4=x_6$ is also a nonsingular toric Fano variety or not. However, there is still a possibility of the decomposition similar to that in Theorem \ref{pikachu} in the case $d\geq 4$.

\begin{Conj}

Let $X$ and $Y$ be nonsingular toric Fano $d$-folds, and $\varphi:Y\longrightarrow X$ an equivariant morphism. Then we have a decomposition of $\varphi$
$$Y=X_r\stackrel{\varphi_r}{\longrightarrow}X_{r-1}\stackrel{\varphi_{r-1}}{\longrightarrow}\cdots\cdots\stackrel{\varphi_2}{\longrightarrow}X_1\stackrel{\varphi_1}{\longrightarrow}X_0=X$$
where $X_i\ (0\leq i\leq r)$ is a nonsingular toric Fano $d$-fold, and each of $\varphi_j\ (1\leq j\leq r)$ is an equivariant blow-up along a $T_{N}$-invariant irreducible closed subvariety of $X_{j-1}$.

\end{Conj}


\section{Program for the classification of toric Fano varieties}\label{class}

\hspace{5mm} In this section, we give a program for the classification of nonsingular toric Fano varieties. This program can be extended to the case of Gorenstein toric Fano varieties endowed with natural resolution of singularities.

First we consider the classification of nonsingular toric Fano $d$-folds. We define the F-equivalence relation again. Let

$$\mathcal{F}_d:=\left\{ \mbox{nonsingular toric Fano $d$-folds}\right\} /\cong.$$

\begin{Def}

{\rm
$X_1$ and $X_2$ in $\mathcal{F}_d$ are said to be {\em F-equivalent} if there exists a sequence of equivariant blow-ups and blow-downs from $X_1$ to $X_2$ through toric {\em Fano} $d$-folds, namely there exist nonsingular toric Fano $d$-folds $Y_0=X_1,Y_1,\ldots,Y_{2l}=X_2$ together with finite successions $Y_j\rightarrow Y_{j-1}$ and $Y_j\rightarrow Y_{j+1}$, for each odd $1\leq j\leq 2l-1$, of equivariant blow-ups through nonsingular toric Fano $d$-folds. We denote the relation by $X_1\stackrel{\rm F}{\sim}X_2$. Then ``$\stackrel{\rm F}{\sim}$'' is obviously an equivalence relation.
}

\end{Def}

By Proposition \ref{kero}, Theorems \ref{fano}, \ref{primup}, \ref{blowdown} and Corollory \ref{primdown}, to get the classification of nonsingular toric Fano $d$-folds, we have only to solve the following problem.

\begin{Prob}\label{pro2}

Get a complete system of representatives for $(\mathcal{F}_d,\stackrel{\rm F}{\sim})$.

\end{Prob}

For Problem \ref{pro2}, we propose the following conjecture.

\begin{Conj}\label{mainconj2}

Any nonsingular toric Fano $d$-fold is either pseudo-symmetric or F-equivalent to the $d$-dimensional projective space ${\bf P}^{d}$, where a nonsingular toric Fano $d$-fold $T_{N}{\rm emb}(\Sigma)$ is {\em pseudo-symmetric} if there exist two $d$-dimensional strongly convex rational polyhedral cones $\sigma$, $\sigma'\in\Sigma(d)$ such that $\sigma=-\sigma':=\{-x\in N_{\bf R} \,|\, x\in\sigma'\}$.

\end{Conj}

If Conjecture \ref{mainconj2} is true, we can get a complete system of representatives for $(\mathcal{F}_d,\stackrel{\rm F}{\sim})$, since pseudo-symmetric ones are already completely classified as follows:

\begin{Def}\label{oyaji}

{\rm
Let $k\in{\bf Z}_{>0}$, $d=2k$ and $\{e_1,\ldots,e_{d}\}$ a basis of $N$. The $2k$-dimensional {\em del Pezzo variety} $V^{2k}$ is the nonsingular toric Fano $2k$-fold corresponding to the Fano polytope in $N_{\bf R}$ defined by
$$\conv\left( \left\{ e_1,\ldots,e_d,-e_1,\ldots,-e_d,e_1+\cdots+e_d,-(e_1+\cdots+e_d)\right\} \right),$$
while the $2k$-dimensional {\em pseudo del Pezzo variety} ${\widetilde{V}}^{2k}$ is the nonsingular toric Fano $2k$-fold corresponding to the Fano polytope in $N_{\bf R}$ defined by
$$\conv\left( \left\{e_1,\ldots,e_d,-e_1,\ldots,-e_d,e_1+\cdots+e_d\right\} \right).$$
}

\end{Def}

\begin{Rem}\label{mame}

{\rm
In the table of Section \ref{four}, $(117)$ is the 4-dimensional pseudo del Pezzo variety $\widetilde{V}^{4}$, while $(118)$ is the 4-dimensional del Pezzo variety $V^{4}$.
}

\end{Rem}

\begin{Thm}[Ewald \cite{ewald1}, Voskresenskij-Klyachko \cite{voskre1}]\label{geraid}

For any pseudo-symmetric toric Fano variety $X$ there exist numbers $s,m,n\in{\bf Z}_{\geq 0},\ k_1,\ldots,k_m,l_1,\ldots l_n\in{\bf Z}_{>0}$ such that
$$X\cong({\bf P}^{1})^{s}\times V^{2k_1}\times\cdots \times V^{2k_m}\times \widetilde{V}^{2l_1}\times\cdots\times\widetilde{V}^{2l_n}$$
where $V^{2k_i}$ is the $2k_i$-dimensional del Pezzo variety while $\widetilde{V}^{2l_j}$ is the $2l_j$-dimensional psudo del Pezzo variety for $1\leq i\leq m,\ 1\leq j \leq n$.

\end{Thm}

Conjecture \ref{mainconj2} is very difficult to deal with in general. So we consider Conjecture \ref{mainconj2} in some special class of nonsingular toric Fano $d$-folds.

\begin{Thm}\label{nobita}

Let $r,a_1,\ldots,a_r$ in ${\bf Z}_{>0}$ and $a_1+\cdots +a_r=d$. Then we have
$${\bf P}^{a_1}\times\cdots\times{\bf P}^{a_r}\stackrel{\rm F}{\sim}{\bf P}^{d}.$$

\end{Thm}

\proof
We prove this by induction on $d$.

Let $\Sigma$ be a fan in $N$ corresponding to the $d$-dimensional projective space and $\G(\Sigma)=\left\{ x_1,\ldots,x_{d+1}\right\}$. Then the primitive relation is $$x_1+\cdots+x_{d+1}=0.$$ By the equivariant blow-up along $\left\{ x_1,\ldots,x_{a_1+1}\right\}$ for $1\leq a_1<d$ we get a fan $\Sigma_1$ in $N$ whose primitive relations are
$$x_1+\cdots+x_{a_1+1}=x_{d+2},\ x_{a_1+2}+\cdots+x_{d+2}=0,$$
where $\G(\Sigma_1)=\G(\Sigma)\cup\left\{ x_{d+2}\right\}$. Moreover, by the equivariant blow-up of $\Sigma_1$ along $\left\{ x_1,x_{a_1+2},\ldots,x_{d+1}\right\}$ we get a fan $\Sigma_2$ in $N$ whose primitive relations are
$$x_1+x_{a_1+2}+\cdots+x_{d+1}=x_{d+3},\ x_2+\cdots+x_{a_1+1}+x_{d+3}=0,\ x_{d+2}+x_{d+3}=x_1,$$
$$x_1+\cdots+x_{a_1+1}=x_{d+2},\ x_{a_1+2}+\cdots+x_{d+2}=0,$$
where $\G(\Sigma_2)=\G(\Sigma_1)\cup\left\{ x_{d+3}\right\}$. Then $T_N{\rm emb}(\Sigma_1)$ and $T_N{\rm emb}(\Sigma_2)$ are nonsingular toric Fano $d$-folds by Theorem \ref{fano}. By Theorem \ref{blowdown} $\Sigma_2$ can be equivariantly blown-down to a fan $\Sigma'$ in $N$ with respect to the primitive relation $x_{d+2}+x_{d+3}=x_1$. The primitive relations of $\Sigma'$ are
$$x_2+\cdots+x_{a_1+1}+x_{d+3}=0,\ x_{a_1+2}+\cdots+x_{d+2}=0,$$
where $\G(\Sigma')=\left\{ x_2,\ldots,x_{d+3}\right\}$. So the toric variety corresponding to $\Sigma'$ is isomorphic to ${\bf P}^{a_1}\times{\bf P}^{d-a_1}$, and we have
$${\bf P}^{d}\stackrel{\rm F}{\sim}{\bf P}^{a_1}\times{\bf P}^{d-a_1}.$$
Then by the induction assumption, we have
$${\bf P}^{d-a_1}\stackrel{\rm F}{\sim}{\bf P}^{a_2}\times\cdots\times{\bf P}^{a_r}.$$
\qed

\bigskip

Next, we consider more complicated nonsingular toric Fano $d$-folds.

\begin{Def}[Batyrev \cite{batyrev3}]

{\rm
Let $\Sigma$ be a finite complete nonsingular fan in $N$. Then $\Sigma$ is called a {\em splitting fan} if for any two distinct primitive collections $P_1$ and $P_2$ in $\PC(\Sigma)$, we have $P_1\cap P_2=\emptyset$.
}

\end{Def}

The following is well-known.

\begin{Thm}[Kleinschmidt \cite{klein1}]\label{ciel}

Let $\Sigma$ be a finite complete nonsingular fan in $N$ and $X:=T_N{\rm emb}(\Sigma)$. If the Picard number of $X$ is two or three, then $X$ is projective. Moreover, if the Picard number of $X$ is two, then $\Sigma$ is a splitting fan.

\end{Thm}

The nonsingular toric $d$-folds corresponding to splitting fans are characterized by the following proposition.

\begin{Prop}[Batyrev \cite{batyrev3}]\label{breeze}

Let $\Sigma$ be a finite complete nonsingular fan in $N$. Then $\Sigma$ is a splitting fan if and only if there exist toric manifolds $X_0,\ldots,X_r$ such that $X_0$ is a projective space, $X_r=T_N{\rm emb}(\Sigma)$ and for $1\leq i\leq r$, $X_i$ is an equivariant projective space bundle over $X_{i-1}$.

\end{Prop}

For any splitting fan $\Sigma$ in $N$, $T_N{\rm emb}(\Sigma)$ is projective by Proposition \ref{breeze}. So the assumption in the following is satisfied.

\begin{Lem}[Batyrev \cite{batyrev3}]

Let $\Sigma$ be a finite complete nonsingular fan in $N$ such that $T_N{\rm emb}(\Sigma)$ is projective. Then there exists a primitive collection $P$ in $\PC(\Sigma)$ such that $\sigma(P)=0$.

\end{Lem}

\begin{Thm}\label{doraemon}

Let $\Sigma$ be a splitting fan in $N$ and let $P=\left\{ x_1,\ldots,x_r\right\}$ be a primitive collection such that $\sigma(P)=0$. If for any primitive collection $P'$ in $\PC(\Sigma)$ such that $\sigma(P')\cap P\neq\emptyset$, there exists $y$ in $P'$ such that $y$ is not in $\sigma(P'')$ for any $P''$ in $\PC(\Sigma)$, then there exists a nonsingular toric Fano $(d-r+1)$-fold $X'$ in ${\mathcal{F}}^{d-r+1}$ such that
$$T_N{\rm emb}(\Sigma)\stackrel{\rm F}{\sim} {\bf P}^{r-1}\times X'.$$

\end{Thm}

\proof
If $\sigma(P')\cap P=\emptyset$ for any primitive collection $P'$ in $\PC(\Sigma)$, then $T_N{\rm emb}(\Sigma)$ is isomorphic to the product as in the statement.

So let $P'=\left\{ y_1,\ldots,y_s\right\}$ be a primitive collection such that $\sigma(P')\cap P\neq\emptyset$ and $x_i$ in $\sigma(P')\cap P$. Then by assumption, there exists $y_j$ in $P'$ such that $y_j$ is not in $\sigma(P'')$ for any $P''$ in $\PC(\Sigma)$. The primitive relations of $\Sigma$ are
$$x_1+\cdots+x_r=0,\ y_1+\cdots+y_s=ax_i+\cdots\ (a>0),\ \ldots.$$
By the equivariant blow-up along $\left\{ x_1,\ldots,x_{i-1},x_{i+1},\ldots,x_r,y_j\right\}$ we get a fan $\Sigma_1$ in $N$ whose primitive relations are
$$x_1+\cdots+x_{i-1}+x_{i+1}+\cdots+x_r+y_j=z,\ x_i+z=y_j,$$
$$y_1+\cdots+y_{j-1}+y_{j+1}+\cdots+y_s+z=(a-1)x_i+\cdots,$$
$$x_1+\cdots+x_r=0,\ y_1+\cdots+y_s=ax_i+\cdots,\ \ldots$$
where $\G(\Sigma_1)=\G(\Sigma)\cup\left\{ z\right\}$ and the first three primitive relations are new. Then $T_N{\rm emb}(\Sigma_1)$ is a nonsingular toric Fano $d$-fold by Theorem \ref{fano}. By Theorem \ref{blowdown} $\Sigma_1$ can be equivariantly blown-down to a fan $\Sigma'$ in $N$ with respect to the primitive relation $x_i+z=y_j$. The primitive relations of $\Sigma'$ are
$$x_1+\cdots+x_r=0,\ y_1+\cdots+y_{j-1}+y_{j+1}+\cdots+y_s+z=(a-1)x_i+\cdots,\ \ldots$$
where $\G(\Sigma')=\left( \G(\Sigma)\setminus\left\{ y_j\right\} \right) \cup\left\{ z\right\}$. Then $T_N{\rm emb}(\Sigma')$ is also a nonsingular toric Fano $d$-fold by Theorem \ref{fano}, and $\Sigma'$ satisfies the assumption of the statement. So we can replace $\Sigma$ by $\Sigma'$ and carry out this operation again. This operation terminates in finite steps and $T_N{\rm emb}(\Sigma')$ becomes a product as in the statement.\qed

\bigskip

By Theorems \ref{nobita} and \ref{doraemon}, we get the following immediately.

\begin{Cor}\label{jiji}

Let $\Sigma$ be a splitting fan in $N$ and let $T_N{\rm emb}(\Sigma)$ be a nonsingular toric Fano $d$-fold. If the Picard number of $T_N{\rm emb}(\Sigma)$ is not greater than three, then $T_N{\rm emb}(\Sigma)$ is F-equivalent to the $d$-dimensional projective space.

\end{Cor}

Next we consider the classification of Gorenstein toric Fano varieties.

Let $\Delta$ be a reflexive polytope in $N_{\bf R}$. For any $\delta\in\Delta(d-1)$, subdivide $\delta$ as
$$\delta=S_{\delta,1}\cup S_{\delta,2}\cup\cdots\cup S_{\delta,k(\delta)}$$
where $S_{\delta,i}\ \left( 1\leq i\leq k(\delta)\right)$ are $(d-1)$-dimensional simplices such that
$$S_{\delta,i}\cap N=S_{\delta,i}(0)\subset\delta\cap N\ \left( 1\leq i\leq k(\delta)\right).$$
Then we can define a finite complete fan $\widetilde{\Sigma(\Delta)}$ in $N$ by
$$\widetilde{\Sigma(\Delta)}:=\left\{ \sigma(S_{\delta,i})\mbox{ and the faces of }\sigma(S_{\delta,i})\; \left| \; \delta\in\Delta(d-1),\ 1\leq i\leq k(\delta)\right. \right\}\cup\{0\}.$$

\begin{Prop}[Batyrev \cite{batyrev2}]\label{nef}

Let $\Delta$ be a reflexive polytope in $N_{\bf R}$. Then there exists a subdivision of $\Sigma(\Delta)$ as above such that $T_N{\rm emb}\left( \widetilde{\Sigma(\Delta)}\right)$ is a projective toric variety with only Gorenstein terminal quotient singularities. Moreover, the equivariant morphism corresponding to this subdivision $\varphi:T_N{\rm emb}\left( \widetilde{\Sigma(\Delta)}\right)\longrightarrow T_N{\rm emb}(\Sigma(\Delta))$ is crepant.

\end{Prop}

\begin{Rem}\label{puffy}

{\rm
In {\rm Proposition \ref{nef}}, if $T_N{\rm emb}\left( \widetilde{\Sigma(\Delta)}\right)$ is nonsingular, then for any $P\in\PC\left( \widetilde{\Sigma(\Delta)}\right)$, we have $\deg P\geq 0$ because $\conv\left( \G\left( \widetilde{\Sigma(\Delta)}\right) \right) =\Delta$. By {\rm Theorem \ref{fano}}, this means that the anticanonical divisor of $T_N{\rm emb}\left( \widetilde{\Sigma(\Delta)}\right)$ is nef.
}

\end{Rem}

\begin{Def}\label{weakfano}

{\rm
Let $X$ be a nonsingular projective algebraic variety. Then $X$ is called a nonsingular {\em weak Fano} variety if the anticanonical divisor $-K_X$ is nef and big.
}

\end{Def}

By the following proposition, the condition ``big'' is automatic in the case of toric varieties.

\begin{Prop}\label{toricweak}

Let $\Sigma$ be a finite projective nonsingular fan in $N$ and $X:=T_{N}{\rm emb}(\Sigma)$. Then the following are equivalent.

\begin{enumerate}

\item $X$ is a nonsingular toric weak Fano variety.
\item The anticanonical divisor $-K_X$ is nef.
\item For any $P\in\PC(\Sigma)$, we have $\deg P\geq 0$.

\end{enumerate}

\end{Prop}

\proof
The equivalence $(2)\Longleftrightarrow (3)$ follows from Theorem \ref{fano}.

Suppose the anticanonical divisor $-K_X$ is nef. Then $\Delta=\conv\left( \G(\Sigma)\right) $ is a reflexive polytope. So we have $(-K_X)^{d}={\rm vol}_d(\Delta^{\ast})>0$. Therefore $-K_X$ is big.\hfill q.e.d.

\bigskip
For the Gorenstein toric Fano varieties endowed with crepant resolutions of singularities as Proposition \ref{nef}, we can consider instead the nonsingular toric weak Fano varieties by Propositions \ref{nef}, \ref{toricweak} and Remark \ref{puffy}. In this case, we can apply the method for nonsingular toric Fano varieties by Theorem \ref{fano} and Proposition \ref{toricweak}. Especially in the cases of $d=2$ and $d=3$, $T_N{\rm emb}\left( \widetilde{\Sigma(\Delta)}\right)$ is always nonsingular.

We introduce the same concepts for nonsingular toric weak Fano $d$-folds as in the case of nonsingular toric Fano $d$-folds. Let
$$\mathcal{F}^{\rm w}_d:=\left\{ \mbox{nonsingular toric weak Fano $d$-folds}\right\} /\cong.$$
First, we define the concept, flop, for nonsingular projective toric $d$-folds.

\begin{Def}

{\rm
Let $X=T_N{\rm emb}(\Sigma)$ be a nonsingular projective toric $d$-fold and $P$ a primitive collection of $\Sigma$ whose primitive relation is
$$x_1+\cdots+x_l=y_1+\cdots+y_l.$$
If $r(P)$ is contained in an extremal ray of ${\bf NE}(X)$, then we can do the following operation. First we blow-up $X$ along $\{y_1,\ldots,y_l\}$, and we get the toric variety $X'=T_N{\rm emb}(\Sigma')$ and the primitive relation of $\Sigma'$, $x_1+\cdots+x_l=z$, where $\G(\Sigma')=\G(\Sigma)\cup\{z\}$. Next we can blow-down $X'$ with respect to $x_1+\cdots+x_l=z$, and we get the toric variety $X^{+}=T_N{\rm emb}(\Sigma^{+})$ and the primitive relation of $\Sigma^{+}$,
$$y_1+\cdots+y_l=x_1+\cdots+x_l,$$
where $\G(\Sigma^{+})=\G(\Sigma)$.
We call this operation {\em flop}.
}

\end{Def}

\begin{Def}

{\rm
$X_1$ and $X_2$ in $\mathcal{F}^{\rm w}_d$ are said to be {\em weakly-F-equivalent} if there exists a sequence of equivariant blow-ups, blow-downs and flops from $X_1$ to $X_2$ through toric {\em weak Fano} $d$-folds, namely there exist nonsingular toric weak Fano $d$-folds $Y_0=X_1,Y_1,\ldots,Y_{3l}=X_2$ together with finite successions $Y_{3j-2}\rightarrow Y_{3j-3}$ and $Y_{3j-2}\rightarrow Y_{3j-1}$, for each $1\leq j\leq l$, of equivariant blow-ups through nonsingular toric Fano $d$-folds, and finite successions $Y_{3k-1}\leftrightarrow Y_{3k}$, for each $1\leq k\leq l$, of flop through nonsingular toric Fano $d$-folds. We denote the relation by $X_1\stackrel{\rm wF}{\sim}X_2$. Then ``$\stackrel{\rm wF}{\sim}$'' is obviously an equivalence relation.
}

\end{Def}

The conjecture for nonsingular toric weak Fano $d$-folds corresponding to Conjecture \ref{mainconj2} is the following.

\begin{Conj}\label{wmainconj}

Any nonsingular toric weak Fano $d$-fold is weakly-F-equivalent to the $d$-dimensional projective space ${\bf P}^{d}$.

\end{Conj}

\begin{Rem}

{\rm
Since the $4$-dimensional pseudo del Pezzo variety and the $4$-dimensional del Pezzo variety can be equivariantly blown-up to nonsingular toric weak Fano $4$-folds, we exclude the pseudo-symmetric toric Fano varieties from Conjecture \ref{wmainconj}.
}

\end{Rem}

We can easily prove Conjectures \ref{mainconj2} and \ref{wmainconj} for $d=2$.

\begin{Thm}\label{coffee}

Any nonsingular toric Fano surface is F-equivalent to the $2$-dimensional projective space ${\bf P}^{2}$, while any nonsingular toric weak Fano surface is weakly-F-equivalent to ${\bf P}^{2}$. Especially, Conjectures {\rm \ref{mainconj2}} and {\rm \ref{wmainconj}} are true for $d=2$, and we get a new method for the classification of Gorenstein toric Fano surfaces by the above discussion.

\end{Thm}

\proof
We prove Theorem \ref{coffee} in the case of nonsingular toric weak Fano surfaces. We can similarly prove Theorem \ref{coffee} in the case of nonsingular toric Fano surfaces.

By Proposition \ref{keylemma} and Example \ref{youming}, if a nonsingular toric weak Fano surface $X$ is not minimal in the sense of equivariant blow-ups, then $X$ can be equivariantly blown-down to nonsingular toric weak Fano surface. On the other hand, the minimal complete nonsingular toric surfaces in the sense of equivariant blow-ups are ${\bf P}^2$ and ${\bf P}_{{\bf P}^{1}}\left( \mathcal{O}_{{\bf P}^{1}}\oplus \mathcal{O}_{{\bf P}^{1}}(a)\right)$ $(a\geq 0$ and $a\neq 1)$ (See Oda \cite{oda2}). So the minimal nonsingular toric weak Fano surfaces in the sense of equivariant blow-ups are
$${\bf P}^{2},\ {{\bf P}^{1}}\times{{\bf P}^{1}}\ \mbox{and}\ {\bf P}_{{\bf P}^{1}}\left( \mathcal{O}_{{\bf P}^{1}}\oplus \mathcal{O}_{{\bf P}^{1}}(2)\right).$$These are weakly-F-equivalent to the $2$-dimensional projective space ${{\bf P}^{2}}$ by easy calculation.\qed


\section{The classification of nonsingular toric Fano $3$-folds}\label{nfano3}

\hspace{5mm} We devote this section to proving Conjecture \ref{mainconj2} for $d=3$. Throughout this section, we assume $d=3$.

\begin{Thm}\label{larc3}

Any nonsingular toric Fano $3$-fold is F-equivalent to the $3$-dimensional projective space ${\bf P}^{3}$. Especially, Conjecture {\rm \ref{mainconj2}} is true for $d=3$, and we get a new method for the classification of nonsingular toric Fano $3$-folds.

\end{Thm}

To prove Theorem \ref{larc3}, we prove the following lemma. For a toric variety $X$, let $\rho(X)$ be the Picard number of $X$.

\begin{Lem}\label{lion}

Let $X=T_{N}{\rm emb}(\Sigma)$ be a nonsingular toric Fano $3$-fold and $\rho(X)\geq 2$. Then there exists a primitive collection $P$ in $\PC(\Sigma)$ such that $\#P=2$.

\end{Lem}

\proof
Suppose there does not exist a primitive collection $P$ in $\PC(\Sigma)$ such that $\#P=2$. Let $\Delta(\Sigma)$ be the Fano polytope corresponding to $X$. Then the $f$-vector of $\Delta(\Sigma)$ is $$\left( \rho(X)+3,\ (\rho(X)+2)(\rho(X)+3)/2,\ f_3\right)$$ by assumption. By the Dehn-Sommerville equalities (see Oda \cite{oda2}), we have $\rho(X)=1$ and $f_3=4$.\qed

\begin{Rem}

{\rm
The method in the proof of Lemma \ref{lion} is not available for $d\geq 4$, because for any $f_0>0$, there always exists a simplicial polytope whose $f$-vector is $$\left( f_0,\ f_0(f_0-1)/2,\ldots\ldots\ \right).$$
}

\end{Rem}

Proof of Theorem \ref{larc3}. \quad
Let $X=T_N{\rm emb}(\Sigma)$ be a nonsingular toric Fano $3$-fold.

If $\rho(X)=2$, then $\Sigma$ is a splitting fan by Theorem \ref{ciel}, and $X$ is F-equivalent to ${\bf P}^{3}$ by Corollary \ref{jiji}.

Suppose $\rho(X)\geq 3$. Then there exists a primitive collection $P$ in $\PC(\Sigma)$ such that $\#P=2$ by Lemma \ref{lion}. By Theorem \ref{fano}, we have two cases.

\begin{enumerate}

\item ``
There exists a primitive collection $P$ in $\PC(\Sigma)$ whose primitive relation is $x_1+x_2=x\ \left( x_1,x_2,x\in\G(\Sigma)\right)$.'' \quad
Because $\deg P=1$, $r(P)$ is contained in an extremal ray of ${\bf NE}(X)$. So $X$ can be equivariantly blown-down with respect to $x_1+x_2=x$. Let $\varphi:X\rightarrow Y$ be the equivariant blow-down with respect to $x_1+x_2=x$. By Proposition \ref{keylemma} and Example \ref{youming}, if $Y$ is not a nonsingular toric Fano $3$-fold, then there exists a primitive collection $P'$ in $\PC(\Sigma)$ whose primitive relation is
$$y_1+y_2=x\quad \left( \{x_1,x_2\}\cap\{y_1,y_2\}=\emptyset\right).$$
By Theorem \ref{blowdown} and the fact $\deg P'=1$, $\{x,x_1,y_1\}$, $\{x,x_1,y_2\}$,$\{x,x_2,y_1\}$ and $\{x,x_2,y_2\}$ generate $3$-dimensional strongly convex rational polyhedral cones of $\Sigma$. Since $\rho(X)\geq 3$, there exists $z$ in $\G(\Sigma)\setminus\{x,x_1,x_2,y_1,y_2\}$. $\{x,z\}$ is obviously a primitive collection of $\Sigma$. If the primitive relation of $\{x,z\}$ is $x+z=z'\ \left( z'\in\G(\Sigma)\right)$, then obviously $X$ can be equivariantly blown-down to a nonsingular toric Fano $3$-fold with respect to $x+z=z'$. If the primitive relation of $\{x,z\}$ is $x+z=0$ and $\rho(X)\geq 4$, then there exists $w$ in $\G(\Sigma)\setminus\{x,x_1,x_2,y_1,y_2,z\}$ and we can replace $z$ by $w$. If the primitive relation of $\{x,z\}$ is $x+z=0$ and $\rho(X)=3$, then the primitive relations of $\Sigma$ are
$$x_1+x_2=x,\ y_1+y_2=x\mbox{ and }x+z=0.$$
Then $\Sigma$ is a splitting fan, and $X$ is F-equivalent to ${\bf P}^{3}$ by Corollary \ref{jiji}.

\item ``For any primitive collection $P$ in $\PC(\Sigma)$ such that $\#P=2$, its primitive relation is $x_1+x_2=0\ \left( x_1,x_2\in\G(\Sigma)\right)$.'' \quad
There exists a primitive relation $x_1+x_2=0$ by Lemma \ref{lion}. Let $\left\{ x_1,x'_1,x''_1\right\}$ generate a $3$-dimensional strongly convex rational polyhedral cone in $\Sigma$, where $x'_1$ and $x''_1$ are in $\G(\Sigma)$. By assumption, there exist $y_1$ and $y_2$ in $\G(\Sigma)\setminus \left\{ x_1,x_2,x'_1,x''_1\right\}$. If $\left\{ y_1,y_2 \right\}$ is not a primitive collection, then $\left\{ x_2,y_1,y_2 \right\}$ generates a $3$-dimensional cone in $\Sigma$. Because $\{x_2,y_1\}$ and $\{x_2,y_2\}$ are also not a primitive collection by assumption, the open set $N_{\bf R}\setminus\left( \cone\left( \{x_2,y_1\}\right) \cup\cone\left( \{x_2,y_2\}\right) \cup\cone\left( \{y_1,y_2\}\right) \right)$ has two connected components. If $\left\{ x_2,y_1,y_2 \right\}$ is a primitive collection, then there exist elements of $\G(\Sigma)$ in both connected components, and there exists a primitive relation like $u_1+u_2=u$. This contradicts the assumption. Either $\left\{ x'_1,y_1\right\}$ or $\left\{ x''_1,y_1\right\}$ is a primitive collection, because otherwise, both $\left\{ x_2,y_1,x'_1 \right\}$ and $\left\{ x_2,y_1,x''_1 \right\}$ generate $3$-dimensional cones in $\Sigma$. So we get primitive relations $y_1+x'_1=0$ and $y_2+x''_1=0$ up to change of the indices. Therefore 
$$\cone\left( \left\{ x_1,x'_1,x''_1 \right\} \right)=-\cone\left( \left\{ x_2,y_1,y_2 \right\} \right),$$
and $T_N{\rm emb}(\Sigma)$ is a pseudo-symmetric toric Fano $3$-fold. Conversely let $\left\{ y_1,y_2 \right\}$ be a primitive collection. Then the corresponding primitive relation is $y_1+y_2=0$ by assumption, and $x_1$, $x_2$, $y_1$ and $y_2$ are contained in a plane. So there exists $z$ in $\G(\Sigma)\setminus \left\{ x_1,x_2,x'_1,x''_1,y_1,y_2\right\}$, and both $\left\{ x'_1,z\right\}$ and $\left\{ x''_1,z\right\}$ are primitive collections. This contradicts the assumption. On the other hand, by Theorem \ref{geraid}, the psudo-symmetric toric Fano $3$-folds are
$${\bf P}^{1}\times{\bf P}^{1}\times{\bf P}^{1},\ {\bf P}^{1}\times V^{2}\mbox{ and }{\bf P}^{1}\times {\widetilde{V}}^{2}.$$
By Definition \ref{oyaji} and Theorem \ref{nobita}, these are F-equivalent to ${\bf P}^{3}$.\qed

\end{enumerate}


\section{The classification of nonsingular toric Fano $4$-folds}\label{nfano4}

\hspace{5mm} In this section, we prove Conjecture \ref{mainconj2} for $d=4$. As a result, we get a new method for the classification of nonsingular toric Fano $4$-folds. Using this method for the classification, we can get the 124 nonsingular toric Fano $4$-folds.

\begin{Thm}\label{larc4}

Any nonsingular toric Fano $4$-fold other than the $4$-dimensional del Pezzo variety $V^{4}$ and the $4$-dimensional pseudo del Pezzo variety ${\widetilde{V}}^{4}$ is F-equivalent to the $4$-dimensional projective space ${\bf P}^{4}$. Especially, Conjecture {\rm \ref{mainconj2}} is true for $d=4$, and we get a new method for the classification of nonsingular toric Fano $4$-folds.

\end{Thm}

We devote the rest of this section to proving Theorem \ref{larc4}. So let $X=T_N{\rm emb}(\Sigma)$ be a nonsingular toric Fano $4$-fold and $\rho=\rho(X)$ the Picard number of $X$.

If $\rho(X)=2$, then $\Sigma$ is a splitting fan by Theorem \ref{ciel}, and $X$ is F-equivalent to ${\bf P}^{4}$ by Corollary \ref{jiji}.

The following theorem holds for nonsingular projective toric $d$-folds for any $d$ whose Picard number is three.

\begin{Thm}[Batyrev \cite{batyrev3}]\label{tmjtmj}

Let $X=T_N{\rm emb}(\Sigma)$ be a nonsingular projective toric $d$-fold such that the Picard number of $X$ is three. Then one of the following holds.

\begin{enumerate}

\item $\Sigma$ is a splitting fan.
\item $\#\PC(\Sigma)=5$.

\end{enumerate}

Moreover, in the case of $(2)$, there exists $(p_0,p_1,p_2,p_3,p_4)\in({\bf Z}_{>0})^{5}$ such that the primitive relation of $\Sigma$ are
$$v_1+\cdots+v_{p_0}+y_{1}+\cdots+y_{p_1}=c_2z_2+\cdots+c_{p_2}z_{p_2}+(b_1+1)t_1+\cdots+(b_{p_3}+1)t_{p_3},$$
$$y_1+\cdots+y_{p_1}+z_1+\cdots+z_{p_2}=u_1+\cdots+u_{p_4},\ z_1+\cdots+z_{p_2}+t_1+\cdots+t_{p_3}=0,$$
$$t_1+\cdots+t_{p_3}+u_1+\cdots+u_{p_4}=y_{1}+\cdots+y_{p_1}\mbox{ and}$$
$$u_1+\cdots+u_{p_4}+v_1+\cdots+v_{p_0}=c_2z_2+\cdots+c_{p_2}z_{p_2}+b_1t_1+\cdots+b_{p_3}t_{p_3},$$
where
$$\G(\Sigma)=\{v_1,\ldots,v_{p_0},y_{1},\ldots,y_{p_1},z_1,\ldots,z_{p_2},t_1,\ldots,t_{p_3},u_1,\cdots,u_{p_4}\},$$
and $c_2,\ldots,c_{p_2},b_1,\ldots,b_{p_3}\in{\bf Z}_{>0}.$

\end{Thm}

For a nonsingular toric Fano variety of any dimension $d$, the following proposition holds.

\begin{Prop}\label{blmmp}

In Theorem {\rm \ref{tmjtmj}}, if X is a nonsingular toric Fano $d$-fold, and $p_1=1$ or $p_4=1$, then $X$ can be equivariantly blown-down to a nonsingular toric Fano $d$-fold.

\end{Prop}

\proof
We prove Proposition \ref{blmmp} in the case of $p_1=1$. We can prove the case of $p_4=1$ similarly.

By assumption, we have the primitive relation of $\Sigma$,
$$t_1+\cdots+t_{p_3}+u_1+\cdots+u_{p_4}=y_{1}.$$
The primitive collections which have a common elements with $\{t_1,\ldots,t_{p3},u_1,\ldots,u_{p_4}\}$ are
$$\{z_1,\ldots,z_{p_2},t_1,\ldots,t_{p_3}\}\mbox{ and }\{u_1,\cdots,u_{p_4},v_1,\ldots,v_{p_0}\}.$$
Since $\{z_1,\ldots,z_{p_2},y_1\}$ and $\{v_1,\ldots,v_{p_0},y_1\}$ are in $\PC(\Sigma)$, by Theorem \ref{blowdown}, $X$ can be equivariantly blown-down to a toric variety $X'$. $X'$ is obviously a nonsingular toric Fano variety by Proposition \ref{keylemma}.\qed

\bigskip

Let $\rho=3$. Since $\#\G(\Sigma)=7$, we get $(p_0,p_1,p_2,p_3,p_4)=(1,1,1,1,3)$, $(1,1,1,2,2)$ or their permutations. By Proposition \ref{blmmp}, if $(p_0,p_1,p_2,p_3,p_4)\neq(1,2,1,1,2)$, then $X$ can be equivariantly blown-down to a nonsingular toric Fano $4$-fold. Let $(p_0,p_1,p_2,p_3,p_4)=(1,2,1,1,2)$. Then the primitive relations of $\Sigma$ are
$$v_1+y_1+y_2=(b_1+1)t_1,\ y_1+y_2+z_1=u_1+u_2,\ z_1+t_1=0,$$
$$t_1+u_1+u_2=y_1+y_2\mbox{ and }u_1+u_2+v_1=b_1t_1,$$
where $b_1=0$ or $1$. If $b_1=0$, then $X$ can be equivariantly blown-down to a nonsingular toric Fano $4$-fold with respect to $v_1+y_1+y_2=t_1$ by Theorem \ref{blowdown} and Proposition \ref{keylemma}. On the other hand, if $b_1=1$, we can show easily that $X$ is F-equivalent to ${\bf P}^{4}$ (see $G_1$ in the table of Section \ref{four}).

Next we consider the case of $\rho\geq4$. We need the following proposition.

\begin{Prop}\label{step1}

Let $X=T_{N}{\rm emb}(\Sigma)$ be a nonsingular toric Fano $4$-fold and $\rho(X)\geq 3$. Then there exists a primitive collection $P$ in $\PC(\Sigma)$ such that $\#P=2$.

\end{Prop}

To prove Proposition \ref{step1}, we have to prove the following three lemmas.

\begin{Lem}\label{tanu1}

Let $X=T_{N}{\rm emb}(\Sigma)$ be a nonsingular toric Fano $4$-fold and $\rho(X)\geq 3$. If there does not exist a primitive collection $P$ in $\PC(\Sigma)$ such that $\#P=2$, then there does not exist a primitive collection $P$ in $\PC(\Sigma)$ such that $\#P=4$.

\end{Lem}

\proof
Suppose there exists a primitive collection $P=\{x_1,x_2,x_3,x_4\}$ in $\PC(\Sigma)$. Then the open set
$$N_{\bf R}\setminus \left( \cone\left( \{x_2,x_3,x_4\}\right) \cup\cone\left( \{x_1,x_3,x_4\}\right) \cup\cone\left( \{x_1,x_2,x_4\}\right) \cup\cone\left( \{x_1,x_2,x_3\}\right) \right)$$
has two connected components. Therefore, since there are two other elements by the assumption $\rho(X)\geq 3$, there exists a primitive relation $P'$ in $\PC(\Sigma)$ such that $\#P'=2$. This contradicts the assumption.\qed

\begin{Lem}\label{tanu2}

Let $X=T_{N}{\rm emb}(\Sigma)$ be a nonsingular toric Fano $4$-fold and $\rho(X)\geq 3$. If there does not exist a primitive collection $P$ in $\PC(\Sigma)$ such that $\#P=2$, then there does not exist a primitive relation of $\Sigma$ of the form
$$x_1+x_2+x_3=ax_4\quad (a=1,2).$$

\end{Lem}

\proof
Suppose there exists a primitive collection $P$ in $\PC(\Sigma)$ whose primitive relation is $x_1+x_2+x_3=ax_4\ (a=1,2)$. If $r(P)$ is contained in an extremal ray of ${\bf NE}(X)$, then there exist $z_1,z_2\in\G(\Sigma)\setminus\{x_1,x_2,x_3,x_4\}$ such that $\{x_i,x_j,x_4,z_k\}$ generate $4$-dimensional strongly convex rational polyhedral cones in $\Sigma$ for $1\leq i<j\leq 3,\ 1\leq k\leq 2$. Since $\#\G(\Sigma)=\rho+4\geq 7$, there exists $w\in\G(\Sigma)\setminus\{x_1,x_2,x_3,x_4,z_1,z_2\}$, and $\{x_4,w\}$ is a primitive collection of $\Sigma$. This contradicts the assumption. So because $\deg P=1$, there does not exists a primitive relation $x_1+x_2+x_3=2x_4$. On the other hand, the primitive relation $x_1+x_2+x_3=x_4$ is represented as the sum of two primitive relations of degree one. By Lemma \ref{tanu1} and assumption, for any primitive collection $P'$ such that $\deg P'=1$, its primitive relation is $y_1+y_2+y_3=y_4+y_5$. Therefore, there exist two primitive relations
$$t_1+t_2+x_1=x_4+s\mbox{ and }s+x_2+x_3=t_1+t_2$$
such that
$$\{t_1,t_2,x_4,s\},\ \{t_1,x_1,x_4,s\},\ \{t_2,x_1,x_4,s\},\ \{s,x_2,t_1,t_2\},\ \{s,x_3,t_1,t_2\}\mbox{ and }\{x_2,x_3,t_1,t_2\}$$
generate $4$-dimensional strongly convex rational polyhedral cones in $\Sigma$, respectively. This is a contradiction because there exist three $4$-dimensional strongly convex rational polyhedral cones generated by $\{t_1,t_2,x_4,s\}$, $\{s,x_2,t_1,t_2\}$ and $\{s,x_3,t_1,t_2\}$, and they contain the $3$-dimensional strongly convex rational polyhedral cone generated by $\{t_1,t_2,s\}$.\qed

\begin{Lem}\label{tanu3}

Let $X=T_{N}{\rm emb}(\Sigma)$ be a nonsingular toric Fano $4$-fold and $\rho(X)\geq 3$. If there does not exist a primitive collection $P$ in $\PC(\Sigma)$ such that $\#P=2$, then there exists a primitive collection $P'=\{x_1,x_2,x_3\}$ in $\PC(\Sigma)$ such that $x_1+x_2+x_3\neq0.$

\end{Lem}

\proof
Suppose $x_1+x_2+x_3=0$ for any primitive collection $P'=\{x_1,x_2,x_3\}$ in $\PC(\Sigma)$. By Lemmas \ref{tanu1} and \ref{tanu2} and assumption, for any primitive collection $P$ in $\PC(\Sigma)$, we have $\#P=3$. If $\Sigma$ is a splitting fan, then $X$ is isomorphic to ${\bf P}^{2}\times{\bf P}^{2}$, and $\rho(X)=2$. So there exist two primitive collections $P_1,P_2$ in $\PC(\Sigma)$ such that $P_1\cap P_2=\emptyset$. If $P_1=\{x_1,x_2,x_3\}$ and $P_2=\{x_1,x_4,x_5\}$, namely $\#(P_1\cap P_2)=1$, then we have $x_2+x_3=x_4+x_5$, and $\{x_2,x_3\}$ or $\{x_4,x_5\}$ in $\PC(\Sigma)$. This contradicts the assumption. The case $P_1=\{x_1,x_2,x_3\}$ and $P_2=\{x_1,x_2,x_4\}$, namely $\#(P_1\cap P_2)=2$, is impossible, because $x_3=x_4$.\qed

\bigskip

Proof of Proposition \ref{step1}. \quad
By Lemmas \ref{tanu2} and \ref{tanu3}, there exists a primitive collection $P$ in $\PC(\Sigma)$ whose primitive relation is $x_1+x_2+x_3=x_4+x_5$. Since $\deg P=1$, we have three $4$-dimensional strongly convex rational polyhedral cones generated by $\{x_i,x_j,x_4,x_5\}$ where $1\leq i<j\leq 3$. There exist $y_1$ and $y_2$ in $\G(\Sigma)\setminus\{x_1,x_2,x_3,x_4,x_5\}$ by the assumption $\rho\geq 3$, and we have $\{x_4,x_5,y_1\}$ and $\{x_4,x_5,y_2\}$ in $\PC(\Sigma)$. If $y_1+x_4+x_5=0$, then $y_2+x_4+x_5\neq0$. Therefore $y_2+x_4+x_5=x_1+x_2$ up to change of indices. This is a contradiction because we have $x_3+y_2=0$ and $\{x_3,y_2\}$ is in $\PC(\Sigma)$. The case $y_1+x_4+x_5\neq0$ is similar.\qed

\bigskip

Let $\rho\geq 4$. Then there exists a primitive collection of $\Sigma$ whose cardinality is two by Proposition \ref{step1}. We divide the proof of Theorem \ref{larc4} for $\rho\geq 4$ into two cases.

\bigskip

$(1)$ ``There exists a primitive relation of $\Sigma$, $x_1+x_2=x$ where $x_1,x_2,x\in\G(\Sigma)$.''

Let $\varphi:X\rightarrow X'$ be the equivariant blow-down with respect to $x_1+x_2=x$. If $X'$ is not a nonsingular toric Fano $4$-fold, then by Proposition \ref{keylemma} and Example \ref{youming}, there exist one of the following primitive relations of $\Sigma$.
$$y_1+y_2+y_3=2x,\ y_1+y_2+y_3=x+x_1\mbox{ and }y_1+y_2=x,$$
where $y_1,y_2,y_3$ in $\G(\Sigma)$.

$(1.1)$ ``$y_1+y_2+y_3=2x$ or $y_1+y_2+y_3=x+x_1$.''\quad Since the degree is one, we have six $4$-dimensional strongly convex rational polyhedral cones generated by $\{x_i,x,y_j,y_k\}$ where $1\leq i\leq 2$, $1\leq j<k\leq 3$. There exist $z_1$ and $z_2$ in $\G(\Sigma)\setminus\{x_1,x_2,x,y_1,y_2,y_3\}$, because $\#\G(\Sigma)=\rho+4\geq 8$, and we have $\{x,z_1\}$ and $\{x,z_2\}$ in $\PC(\Sigma)$. Therefore we get a primitive relation of $\Sigma$ of the form
$$x+z_1=w\quad \left( w\in\{x_1,x_2,y_1,y_2,y_3\}\right)$$
up to change of indices. Let $\varphi:X\rightarrow X''$ be the equivariant blow-down with respect to $x+z_1=w$.

$(1.1.1)$ ``$w\neq x_1$ or $y_1+y_2+y_3=2x$ is a primitive relation of $\Sigma$.''\quad Then $X''$ is obviously a nonsingular toric Fano $4$-fold.

$(1.1.2)$ ``$w=x_1$ and $y_1+y_2+y_3=x+x_1$ is a primitive relation of $\Sigma$.''\quad In this case, $X''$ is not a nonsingular toric Fano 4-fold by Proposition \ref{keylemma} and Example \ref{youming}. Since $\rho\geq 4$, there exists $t\in \G(\Sigma)\setminus \{ x_1,x_2,x,y_1,y_2,y_3,z_1\}$. So we get one of the following primitive relations of $\Sigma$,
$$t+x_1=y_1,\ t+x_1=x_2\mbox{ and }t+x_1=z_1,$$
up to change of indices. Let $\varphi':X\rightarrow X'''$ be the equivariant blow-down with respect to this primitive relation. Then $X'''$ is obviously a nonsingular toric Fano $4$-fold.

$(1.2)$ ``$y_1+y_2=x$''\quad Since the degree is one, there exist two elements $z_1$ and $z_2$ in $\G(\Sigma)\setminus\{x_1,x_2,x,y_1,y_2\}$, and we have eight $4$-dimensional strongly convex rational polyhedral cones generated by $\{x_i,x,y_j,z_k\}$ where $1\leq i,j,k \leq 2$. There exist $w$ in $\G(\Sigma)\setminus\{x_1,x_2,x,y_1,y_2,z_1,z_2\}$, because $\#\G(\Sigma)=\rho+4\geq 8$, and $P=\{x,w\}$ is a primitive collection of $\Sigma$.

$(1.2.1)$ ``The primitive relation of $P$ is $x+w=t$ where $t$ in $\{z_1,z_2\}$.''\quad Let $\varphi:X\rightarrow X''$ be the equivariant blow-down with respect to $x+w=t$. Then $X''$ is obviously a nonsingular toric Fano $4$-fold.

$(1.2.2)$ ``The primitive relation of $P$ is $x+w=t$ where $t$ in $\{x_1,x_2,y_1,y_2\}$.''\quad Let $\varphi:X\rightarrow X''$ be the equivariant blow-down with respect to $x+w=t$. If $X''$ is not a nonsingular toric Fano $4$-fold, then we obviously have a primitive relation $z_1+z_2=t$ by Proposition \ref{keylemma}. We may let $t=x_2$ without loss of generality. Then we have four $4$-dimensional strongly convex rational polyhedral cones generated by $\{ x_2,y_i,y_j,w\}$ where $1\leq i,j\leq 2$. $\{x_1,x,y_1,z_1\}$ is a ${\bf Z}$-basis of $N$, and using this basis, we have
$$x_2=-x_1+x,\ y_2=x-y_1,\ z_2=-x_1+x-z_1\mbox{ and }w=-x_1.$$
Since the coefficient of $x$ in none of these relation is negative, there exist $u$ in $\G(\Sigma)\setminus\{x_1,x_2,x,y_1,y_2,z_1,z_2,w\}$ by the completeness of $\Sigma$, and we have two primitive collections $\{x,u\}$ and $\{x_2,u\}$ in $\PC(\Sigma)$. Therefore, we have a primitive relation, $x+u=s$ or $x_2+u=s$ where $s$ is in $\{x_1,y_1,y_2,z_1,z_2,w\}$. Let $\varphi:X\rightarrow X''$ be the equivariant blow-down with respect to $x+u=s$. Then $X''$ is obviously a nonsingular toric Fano $4$-fold. The case of the blow-down with respect to $x_2+u=s$ is similar.

$(1.2.3)$ ``The primitive relation of $P$ is $x+w=0$.''\quad If $\rho\geq 5$, then there exist $v$ in $\G(\Sigma)\setminus\{x_1,x_2,x,y_1,y_2,z_1,z_2,w,v\}$, and we have the primitive relation $x+v\neq 0$. In this case, we can use the same method as in $(1.2.1)$ or $(1.2.2)$.

So let $\rho=4$ and $\G(\Sigma)=\{x_1,x_2,x,y_1,y_2,z_1,z_2,w\}$. Then either $\{z_1,z_2\}$ or $\{x,z_1,z_2\}$ is a primitive collection of $\Sigma$.

$(1.2.3.1)$ ``$z_1+z_2=0$ is a primitive relation of $\Sigma$.''\quad $X$ is obviously a nonsingular toric Fano $4$-fold in this case. The primitive relations of $\Sigma$ are
$$x_1+x_2=x,\ y_1+y_2=x,\ x+w=0\mbox{ and }z_1+z_2=0.$$
$\Sigma$ is a splitting fan, and $X$ is F-equivalent to ${\bf P}^{4}$ by Theorems \ref{nobita} and \ref{doraemon}.

$(1.2.3.2)$ ``$z_1+z_2=x$ is a primitive relation of $\Sigma$.''\quad $X$ is obviously a nonsingular toric Fano $4$-fold in this case. The primitive relations of $\Sigma$ are
$$x_1+x_2=x,\ y_1+y_2=x,\ x+w=0\mbox{ and }z_1+z_2=x.$$
$\Sigma$ is a splitting fan, and $X$ is F-equivalent to ${\bf P}^{4}$ by Theorems \ref{nobita} and \ref{doraemon}.

$(1.2.3.3)$ ``$z_1+z_2=t$ is a primitive relation of $\Sigma$, where $t$ in $\{x_1,x_2,y_1,y_2,w\}$.''\quad Let $\varphi:X\rightarrow X''$ be the equivariant blow-down with respect to $z_1+z_2=t$. Then $X''$ is obviously a nonsingular toric Fano $4$-fold by Proposition \ref{keylemma} and Example \ref{youming}.

$(1.2.3.4)$ ``$z_1+z_2+x=0$ is a primitive relation of $\Sigma$.''\quad This is impossible, because $z_1+z_2=-x=w$, and $\{z_1,z_2\}$ is a primitive collection of $\Sigma$.

$(1.2.3.5)$ ``$z_1+z_2+x=ax_1$ is a primitive relation of $\Sigma$, where $a=1$ or $2$.''\quad Since $ax_1+w=z_1+z_2$, $\{t,w\}$ is a primitive collection of $\Sigma$. There exists $u$ in $\{x_2,y_1,y_2,z_1,z_2\}$ such that the primitive relation of $\{x_1,w\}$ is $x_1+w=u$, because $x+w=0$. Since $x_1-x-u=0$, we have $u=x_2$. Because otherwise, $\{x_1,x,u\}$ is a part of a ${\bf Z}$-basis of $N$. However, this contradicts the fact $x_1+x_2=x$. We can replace $x_1$ by $x_2$, $y_1$ or $y_2$, and repeat the same argument.

$(1.2.3.6)$ ``$z_1+z_2+x=aw$ is a primitive relation of $\Sigma$, where $a=1$ or $2$.''\quad We have $z_1+z_2=aw-x=(a+1)w$. This is a contradiction.

$(1.2.3.7)$ ``$z_1+z_2+x=x_i+y_j$ is a primitive relation of $\Sigma$, where $1\leq i,j\leq 2$.''\quad $X$ is obviously a nonsingular toric Fano $4$-fold. We can show easily that $X$ is F-equivalent to ${\bf P}^{4}$ (See $M_2$ in the table of Section \ref{four}).

\bigskip

$(2)$ ``There does not exist a primitive collection $P=\{x_1,x_2\}$ in $\PC(\Sigma)$ whose primitive relation is $x_1+x_2\neq0$.''

We need the following lemma. This lemma can be proved in the same way as Lemmas \ref{tanu1} and \ref{tanu2}.

\begin{Lem}\label{nanako}

Let $X=T_{N}{\rm emb}(\Sigma)$ be a nonsingular toric Fano $4$-fold and $\rho(X)\geq 4$. If there does not exist a primitive collection $P=\{x_1,x_2\}$ in $\PC(\Sigma)$ whose primitive relation is $x_1+x_2\neq0$, then the following hold.

\begin{enumerate}

\item There does not exist a primitive collection $P$ in $\PC(\Sigma)$ such that $\#P=4$.
\item There does not exist a primitive relation of $\Sigma$ of the form
$$x_1+x_2+x_3=ax_4\quad (a=1,2).$$
\end{enumerate}

\end{Lem}

Since $\rho\geq 4$, there exists a primitive collection $P=\{x_1,x_2\}$ in $\PC(\Sigma)$ whose primitive relation is $x_1+x_2=0$, by Proposition \ref{step1}. We fix this $P$.

$(2.1)$ ``$r(P)$ is contained in an extremal ray of ${\bf NE}(X)$.''\quad By toric Mori theory, there exists a nonsingular projective toric $3$-fold $Y=T_{N}{\rm emb}(\Sigma^{\ast})$ such that $X$ is an equivariant ${\bf P}^{1}$-bundle over $Y$, $\G(\Sigma^{\ast})\subset\G(\Sigma)$ and if $P^{\ast}$ is a primitive collection of $\Sigma^{\ast}$, then $P^{\ast}$ is also a primitive collection of $\Sigma$. Let $\#\G(\Sigma^{\ast})=n$ and $n_0$ the number of the primitive collections of $\Sigma^{\ast}$ whose cardinality is two. Then the $f$-vector of the $3$-dimensional simplicial convex polytope corresponding to $\Sigma^{\ast}$ is $(n,n(n-1)/2-n_0,f_2)$. By the Dehn-Sommerville equalities (see Oda \cite{oda2}), we have $n_0=(n-3)(n-4)/2$. So by assumption, we have $n_0=(n-3)(n-4)/2\leq n/2$. Since $\rho\geq 4$, we have $n=6$, and the primitive relations of $\Sigma$ are
$$x_1+x_2=0,\ x_3+x_4=0,\ x_5+x_6=0\mbox{ and }x_7+x_8=0.$$
$X$ is ${\bf P}^{1}\times{\bf P}^{1}\times{\bf P}^{1}\times{\bf P}^{1}$ and F-equivalent to ${\bf P}^{4}$ by Theorem \ref{nobita}.

$(2.2)$ ``$r(P)$ is not contained in an extremal ray of ${\bf NE}(X)$.''\quad By Lemma \ref{nanako}, there exist two primitive relations of $\Sigma$ of the form
$$x_1+y_1+y_2=z_1+z_2\mbox{ and }x_2+z_1+z_2=y_1+y_2,$$
with $y_1,y_2,z_1,z_2$ in $\G(\Sigma)$. We have five $4$-dimensional strongly convex rational polyhedral cones of $\Sigma$ generated by
$$\{x_1,y_1,z_1,z_2\},\ \{x_1,y_2,z_1,z_2\},\ \{y_1,y_2,z_1,z_2\},\ \{x_2,y_1,y_2,z_1\}\mbox{ and }\{x_2,y_1,y_1,z_2\}.$$
By the assumption $\rho\geq 4$, there exists $w$ in $\G(\Sigma)\setminus\{x_1,x_2,y_1,y_2,z_1,z_2\}$ such that either $\{z_1,z_2,w\}$ or $\{z_1,z_2,w\}$ is a primitive collection of $\Sigma$, because there exists at most one primitive collection among $\{z_1,w\}$, $\{z_2,w\}$, $\{y_1,w\}$ and $\{y_1,w\}$, and the others generate $2$-dimensional strongly convex rational polyhedral cones of $\Sigma$. If $w+z_1+z_2=0$ is a primitive relation, then we have $y_1+y_2+w=x_2$. So by assumption, $\{y_1,y_2\}$, $\{y_1,w\}$ and $\{y_2,w\}$ are not primitive collections. Therefore, $\{y_1,y_2,w\}$ is a primitive collection of $\Sigma$. This contradicts Lemma \ref{nanako}.

By the above discussion, we have the primitive relations $w+z_1+z_2=t_1+t_2$ and $w+y_1+y_2=s_1+s_2$, where the possibilities for $\{t_1,t_2\}$ are $\{x_1,y_1\}$ and $\{x_1,y_2\}$, while the possibilities for $\{s_1,s_2\}$ are $\{x_2,z_1\}$ and $\{x_2,z_2\}$. So we have $4\leq\rho\leq 6$.

$(2.2.1)$ ``$\rho=4$''\quad $X$ is obviously a nonsingular toric Fano $4$-fold. We can show easily that $X$ is F-equivalent to ${\bf P}^{4}$ (See $M_1$ in the table of Section \ref{four}).

$(2.2.2)$ ``$\rho=5$''\quad $X$ is the $4$-dimensional pseudo del Pezzo variety. Moreover, $X$ is not F-equivalent to ${\bf P}^{4}$ (See $(117)$ in the table of Section \ref{four}). The primitive relations of $\Sigma$ are
$$x_0+x_4=0,\ x_1+x_5=0,\ x_2+x_6=0,\ x_3+x_7=0,$$
$$x_0+x_1+x_2=x_7+x_8,\ x_0+x_1+x_3=x_6+x_8,\ x_0+x_2+x_3=x_5+x_8,\ x_1+x_2+x_3=x_4+x_8,$$
$$x_4+x_5+x_8=x_2+x_3,\ x_4+x_6+x_8=x_1+x_3,\ x_4+x_7+x_8=x_1+x_2,\ x_5+x_6+x_8=x_0+x_3,$$
$$x_5+x_7+x_8=x_0+x_2,\ x_6+x_7+x_8=x_0+x_1,$$
where $\G(\Sigma)=\{ x_0,x_1,x_2,x_3,x_4,x_5,x_6,x_7,x_8\}$.

$(2.2.3)$ ``$\rho=6$''\quad $X$ is the $4$-dimensional del Pezzo variety. Moreover, $X$ is not F-equivalent to ${\bf P}^{4}$ (See $(118)$ in the table of Section \ref{four}). The primitive relations of $\Sigma$ are
$$x_0+x_4=0,\ x_1+x_5=0,\ x_2+x_6=0,\ x_3+x_7=0,\ x_8+x_9=0,$$
$$x_0+x_1+x_2=x_7+x_8,\ x_0+x_1+x_3=x_6+x_8,\ x_0+x_2+x_3=x_5+x_8,\ x_1+x_2+x_3=x_4+x_8,$$
$$x_0+x_1+x_9=x_6+x_7,\ x_0+x_2+x_9=x_5+x_7,\ x_0+x_3+x_9=x_5+x_6,\ x_1+x_2+x_9=x_4+x_7,$$
$$x_1+x_3+x_9=x_4+x_6,\ x_2+x_3+x_9=x_4+x_5,\ x_4+x_5+x_6=x_3+x_9,\ x_4+x_5+x_7=x_2+x_9,$$
$$x_4+x_6+x_7=x_1+x_9,\ x_5+x_6+x_7=x_0+x_9,\ x_4+x_5+x_8=x_2+x_3,\ x_4+x_6+x_8=x_1+x_3,$$
$$x_4+x_7+x_8=x_1+x_2,\ x_5+x_6+x_8=x_0+x_3,\ x_5+x_7+x_8=x_0+x_2,\ x_6+x_7+x_8=x_0+x_1,$$
where $\G(\Sigma)=\{ x_0,x_1,x_2,x_3,x_4,x_5,x_6,x_7,x_8,x_9\}$.


\section{Equivariant blow-up relations among nonsingular toric Fano $4$-folds}\label{four}       

\hspace{5mm} In this section, we give all the equivariant blow-up relations among nonsingular toric Fano $4$-folds using the results of Sections \ref{prim}, \ref{blow}, \ref{class} and \ref{nfano4}.
\ In Table $1$, we use the same notation as in Batyrev \cite{batyrev4}, and $i$-blow-up means the equivariant blow-up along a $T_{N}$-invariant irreducible closed subvariety of codimension $i$.
\setlongtables
\bigskip

\bigskip

\begin{center}
Table 1: equivariant blow-up relations among nonsingular toric Fano $4$-folds
\end{center}

\begin{longtable}{|c||l||c|}

\hline
\endhead
\hline
\endfoot

\hspace{1cm} &  equivariant blow-up\hspace{6.7cm} & notation\hspace{0.5cm} \\ \hline
(1) & none & ${\bf P}^{4}$ \\ \hline
(2) & none & $B_1$ \\ \hline
(3) & none & $B_2$ \\ \hline
(4) & 4-blow-up of ${\bf P}^{4}$ & $B_3$ \\ \hline
(5) & none & $B_4$ \\ \hline
(6) & 2-blow-up of ${\bf P}^{4}$ & $B_5$ \\ \hline
(7) & none & $C_1$ \\ \hline
(8) & 3-blow-up of ${\bf P}^{4}$ & $C_2$ \\ \hline
(9) & none & $C_3$ \\ \hline
(10) & none & $C_4$ \\ \hline
(11) & 2-blow-up of $B_1,\ B_2$ & $E_1$ \\ \hline
(12) & 2-blow-up of $B_2,\ B_3$ & $E_2$ \\ \hline
(13) & 2-blow-up of $B_3,\ B_4$,\quad 4-blow-up of $B_5$ & $E_3$ \\ \hline
(14) & none & $D_1$ \\ \hline
(15) & 2-blow-up of $C_1$ & $D_2$ \\ \hline
(16) & none & $D_3$ \\ \hline
(17) & 2-blow-up of $B_2$ & $D_4$ \\ \hline
(18) & none & $D_5$ \\ \hline
(19) & 2-blow-up of $C_3$ & $D_6$ \\ \hline
(20) & none & $D_7$ \\ \hline
(21) & 2-blow-up of $C_2$,\quad 3-blow-up of $B_3$ & $D_8$ \\ \hline
(22) & none & $D_9$ \\ \hline
(23) & 2-blow-up of $B_5$,\quad 3-blow-up of $B_3$ & $D_{10}$ \\ \hline
(24) & 2-blow-up of $B_5,\ C_2$ & $D_{11}$ \\ \hline
(25) & 3-blow-up of $B_4$ & $D_{12}$ \\ \hline
(26) & none & $D_{13}$ \\ \hline
(27) & 2-blow-up of $B_4$ & $D_{14}$ \\ \hline
(28) & 2-blow-up of $C_4$ & $D_{15}$ \\ \hline
(29) & 2-blow-up of $C_3$ & $D_{16}$ \\ \hline
(30) & 2-blow-up of $B_5$ & $D_{17}$ \\ \hline
(31) & 2-blow-up of $C_1$ & $D_{18}$ \\ \hline
(32) & 2-blow-up of $C_2$,\quad 3-blow-up of $B_5$ & $D_{19}$ \\ \hline
(33) & none & $G_1$ \\ \hline
(34) & 2-blow-up of $C_2$,\quad 3-blow-up of $C_1$ & $G_2$ \\ \hline
(35) & 3-blow-up of $C_3$ & $G_3$ \\ \hline
(36) & 2-blow-up of $C_2$,\quad 3-blow-up of $C_3$ & $G_4$ \\ \hline
(37) & 2-blow-up of $C_3$,\quad 3-blow-up of $C_4$ & $G_5$ \\ \hline
(38) & 2-blow-up of $C_4$ & $G_6$ \\ \hline
(39) & 2-blow-up of $D_2$ & $H_1$ \\ \hline
(40) & 2-blow-up of $D_3$ & $H_2$ \\ \hline
(41) & 2-blow-up of $D_1,\ D_5$ & $H_3$ \\ \hline
(42) & 2-blow-up of $D_8,\ D_9$ & $H_4$ \\ \hline
(43) & 2-blow-up of $D_6,\ D_{12},\ D_{16}$ & $H_5$ \\ \hline
(44) & 2-blow-up of $D_3,\ D_9$ & $H_6$ \\ \hline
(45) & 2-blow-up of $D_2,\ D_5,\ D_{18}$ & $H_7$ \\ \hline
(46) & 2-blow-up of $D_{13},\ D_{15}$ & $H_8$ \\ \hline
(47) & 2-blow-up of $D_8,\ D_{12},\ D_{19}$,\quad 3-blow-up of $E_3$ & $H_9$ \\ \hline
(48) & 2-blow-up of $D_9,\ D_{16}$ & $H_{10}$ \\ \hline
(49) & none & $L_1$ \\ \hline
(50) & 2-blow-up of $D_7$ & $L_2$ \\ \hline
(51) & 2-blow-up of $D_6$ & $L_3$ \\ \hline
(52) & 2-blow-up of $D_8,\ D_{10},\ D_{11}$ & $L_4$ \\ \hline
(53) & none & $L_5$ \\ \hline
(54) & 2-blow-up of $D_{12},\ D_{14}$ & $L_6$ \\ \hline
(55) & 2-blow-up of $D_{15}$ & $L_7$ \\ \hline
(56) & none & $L_8$ \\ \hline
(57) & 2-blow-up of $D_{13}$ & $L_9$ \\ \hline
(58) & 2-blow-up of $D_{10},\ D_{17}$ & $L_{10}$ \\ \hline
(59) & 2-blow-up of $D_{14}$ & $L_{11}$ \\ \hline
(60) & 2-blow-up of $D_{11},\ D_{17},\ D_{19}$ & $L_{12}$ \\ \hline
(61) & 2-blow-up of $D_7$ & $L_{13}$ \\ \hline
(62) & 2-blow-up of $D_4$ & $I_1$ \\ \hline
(63) & 2-blow-up of $D_1,\ D_6$ & $I_2$ \\ \hline
(64) & 2-blow-up of $D_3,\ D_8$ & $I_3$ \\ \hline
(65) & 2-blow-up of $D_{10}$ & $I_4$ \\ \hline
(66) & 2-blow-up of $E_2,\ D_4,\ D_{10}$ & $I_5$ \\ \hline
(67) & 2-blow-up of $D_{10}$,\quad 3-blow-up of $D_{11}$ & $I_6$ \\ \hline
(68) & 2-blow-up of $D_5,\ D_{12}$ & $I_7$ \\ \hline
(69) & 2-blow-up of $D_8,\ D_{16}$,\ $G_4$ & $I_8$ \\ \hline
(70) & 2-blow-up of $D_{14}$,\quad 3-blow-up of $D_7$ & $I_9$ \\ \hline
(71) & 2-blow-up of $D_6,\ D_{15},\ G_5$ & $I_{10}$ \\ \hline
(72) & 2-blow-up of $D_9,\ D_{12}$ & $I_{11}$ \\ \hline
(73) & 2-blow-up of $D_{15},\ D_{19},\ G_6$,\quad 3-blow-up of $D_{11}$ & $I_{12}$ \\ \hline
(74) & 2-blow-up of $D_{12},\ D_{13}$,\quad 3-blow-up of $D_{14}$ & $I_{13}$ \\ \hline
(75) & 2-blow-up of $E_3,\ D_{10},\ D_{14}$,\quad 3-blow-up of $D_{17}$ & $I_{14}$ \\ \hline
(76) & 2-blow-up of $D_{18},\ D_{19},\ G_2$ & $I_{15}$ \\ \hline
(77) & none & $M_1$ \\ \hline
(78) & none & $M_2$ \\ \hline
(79) & 2-blow-up of $G_3,\ G_5$ & $M_3$ \\ \hline
(80) & 2-blow-up of $G_3$ & $M_4$ \\ \hline
(81) & 2-blow-up of $G_4,\ G_6$ & $M_5$ \\ \hline
(82) & 2-blow-up of $G_1,\ G_3$ & $J_1$ \\ \hline
(83) & 2-blow-up of $G_3$,\quad 3-blow-up of $G_5$ & $J_2$ \\ \hline
(84) & 2-blow-up of $L_2$ & $Q_1$ \\ \hline
(85) & 2-blow-up of $H_4,\ L_4$ & $Q_2$ \\ \hline
(86) & 2-blow-up of $L_1,\ L_5$ & $Q_3$ \\ \hline
(87) & 2-blow-up of $L_3$ & $Q_4$ \\ \hline
(88) & 2-blow-up of $H_5,\ L_3,\ L_6$ & $Q_5$ \\ \hline
(89) & 2-blow-up of $L_6$ & $Q_6$ \\ \hline
(90) & 2-blow-up of $L_7$ & $Q_7$ \\ \hline
(91) & 2-blow-up of $L_5,\ L_9$ & $Q_8$ \\ \hline
(92) & 2-blow-up of $L_3,\ L_7,\ I_{10}$ & $Q_9$ \\ \hline
(93) & 2-blow-up of $H_8,\ L_7,\ L_9$ & $Q_{10}$ \\ \hline
(94) & 2-blow-up of $L_8,\ L_9$ & $Q_{11}$ \\ \hline
(95) & 2-blow-up of $L_{10},\ L_{12},\ I_6$ & $Q_{12}$ \\ \hline
(96) & 2-blow-up of $L_2,\ L_5,\ L_{13}$ & $Q_{13}$ \\ \hline
(97) & 2-blow-up of $H_9,\ L_4,\ L_6,\ L_{12},\ I_{14}$ & $Q_{14}$ \\ \hline
(98) & 2-blow-up of $L_6,\ L_9,\ L_{11},\ I_{13}$ & $Q_{15}$ \\ \hline
(99) & 2-blow-up of $L_{11},\ L_{13},\ I_9$ & $Q_{16}$ \\ \hline
(100) & 2-blow-up of $L_7,\ L_{12},\ I_{12}$ & $Q_{17}$ \\ \hline
(101) & 2-blow-up of $H_1,\ H_3,\ H_7$ & $K_1$ \\ \hline
(102) & 2-blow-up of $H_2,\ H_6,\ H_{10}$ & $K_2$ \\ \hline
(103) & 2-blow-up of $H_4,\ H_5,\ H_9$ & $K_3$ \\ \hline
(104) & 2-blow-up of $H_8$ & $K_4$ \\ \hline
(105) & 2-blow-up of $M_3$ & $R_1$ \\ \hline
(106) & 2-blow-up of $M_2,\ M_4$ & $R_2$ \\ \hline
(107) & 2-blow-up of $M_1,\ M_4$ & $R_3$ \\ \hline
(108) & 2-blow-up of $I_{11},\ I_{13}$ & \\ \hline
(109) & 2-blow-up of $Q_1,\ Q_3,\ Q_{13}$ & $U_1$ \\ \hline
(110) & 2-blow-up of $Q_2,\ Q_5,\ Q_{14},\ K_3$ & $U_2$ \\ \hline
(111) & 2-blow-up of $Q_4,\ Q_9$ & $U_3$ \\ \hline
(112) & 2-blow-up of $Q_{10},\ K_4$ & $U_4$ \\ \hline
(113) & 2-blow-up of $Q_{11}$ & $U_5$ \\ \hline
(114) & 2-blow-up of $Q_6,\ Q_8,\ Q_{15}$ & $U_6$ \\ \hline
(115) & 2-blow-up of $Q_7,\ Q_{12},\ Q_{17}$ & $U_7$ \\ \hline
(116) & 2-blow-up of $Q_{16}$ & $U_8$ \\ \hline
(117) & none (See Definition \ref{oyaji} and Remark \ref{mame}) & ${\widetilde{V}}^{4}$ \\ \hline
(118) & none (See Definition \ref{oyaji} and Remark \ref{mame}) & $V^{4}$ \\ \hline
(119) & 2-blow-up of $Q_{10},\ Q_{11}$ & $S_2\times S_2$ \\ \hline
(120) & 2-blow-up of $U_4,\ U_5,\ S_2\times S_2$ & $S_2\times S_3$ \\ \hline
(121) & 2-blow-up of $S_2\times S_3$ & $S_3\times S_3$ \\ \hline
(122) & 2-blow-up of $G_6$ & $Z_1$ \\ \hline
(123) & 2-blow-up of $G_4$ & $Z_2$ \\ \hline
(124) & 2-blow-up of $Z_1$ (See Example \ref{newfano}) & $W$ \\ \hline

\end{longtable}

\bigskip

\begin{flushleft}
\begin{sc}
Mathematical Institute \\
Tohoku University \\
Sendai 980-8578 \\
Japan
\end{sc}

\medskip
{\it E-mail address}: 96m16@math.tohoku.ac.jp
\end{flushleft}

\end{document}